\documentclass[12pt]{article}  

\usepackage{amsthm}
\renewcommand{\proof}{ \noindent{\it Proof:\ \ }}
\usepackage{hyperref}
\usepackage{amsmath}
\usepackage{amssymb}
\usepackage{amsfonts}
\usepackage[latin1]{inputenc}
\usepackage{enumerate}

\let\oldmarginpar\marginpar
\renewcommand\marginpar[1]{${}^\clubsuit$\oldmarginpar[\raggedleft\scriptsize\sf #1]{\raggedright\scriptsize\sf #1}}
\numberwithin{equation}{section}
\newtheorem{theorem}[equation]{Theorem}

\newtheorem{proposition}[equation]{Proposition}
\newtheorem{lemma}[equation]{Lemma}
\newtheorem{corollary}[equation]{Corollary}
\newtheorem{remark}[equation]{Remark}

\newtheorem{main}{Theorem}
\newtheorem{cor}{Corollary}

\numberwithin{equation}{section}

\def\qed{\ifhmode\unskip\nobreak\fi\ifmmode\ifinner\else\hskip5 pt
\fi\fi\hbox{\hskip5 pt \vrule width4 pt height6 pt depth1.5 pt
\hskip 1pt }}

\let\oldmarginpar\marginpar
\renewcommand\marginpar[1]{${}^\clubsuit$\oldmarginpar[\raggedleft\scriptsize\sf #1]{\raggedright\scriptsize\sf #1}}
\newcommand{\cref}[1]{Corollary~\ref{#1}}

\renewcommand{\eqref}[1]{(\ref{#1})}
\newcommand{\lref}[1]{Lemma~\ref{#1}}
\newcommand{\pref}[1]{Proposition~\ref{#1}}

\newcommand{\tref}[1]{Theorem~\ref{#1}}

\def\Bbb#1{\mathbb#1}
\newcommand{\R}{\Bbb{R}}

\newcommand{\Z}{\Bbb{Z}}
\newcommand{\N}{\Bbb{N}}
\newcommand{\Sp}{\Bbb{S}}

\newcommand{\C}{\Bbb{C\,}}
\newcommand{\CP}{\Bbb{CP\,}}
\newcommand{\e}{\epsilon}
\newcommand{\QH}{\Bbb{H}}
\newcommand{\HP}{\Bbb{HP\,}}

\newcommand{\fl}{{\mathfrak{l}}}
\newcommand{\fg}{{\mathfrak{g}}}
\newcommand{\fh}{{\mathfrak{h}}}
\newcommand{\fo}{{\mathfrak{o}}}
\newcommand{\fp}{{\mathfrak{p}}}
\newcommand{\fs}{{\mathfrak{s}}}
\newcommand{\ft}{{\mathfrak{t}}}
\newcommand{\fu}{{\mathfrak{u}}}
\newcommand{\fsu}{{\mathfrak{su}}}
\newcommand{\sg}{{\rm sign}}

\newcommand{\ra}{\rangle}
\newcommand{\la}{\langle}
\newcommand{\ml}{\langle}
\newcommand{\mr}{\rangle}

\newcommand{\h}{\mathcal{H}}
\newcommand{\diag}{{\rm diag}}
\newcommand{\fsp}{{\mathfrak{sp}}}
\newcommand{\Pf}{{\rm Pf}}

\newcommand{\im}{{\rm Im \,}}
\newcommand{\rk}{{\rm rank \,}}

\renewcommand{\v}{\mathcal{V}}
\renewcommand{\h}{\mathcal{H}}

\usepackage{fancyhdr}
\fancypagestyle{plain}{
\cfoot{}
\lfoot{
\footnotesize
Both authors were partially supported by CNPq-Brazil and the second
one by a grant from the National Science Foundation.
\begin{center} \thepage \end{center}
}

}

\begin{document}
\title{Topological obstructions to fatness}  
\author{Luis A. Florit and Wolfgang Ziller}

%

\date{}
\maketitle

There are few known examples of compact Riemannian manifolds with
positive sectional curvature; see \cite{zi2} for a survey. All of them,
apart from some rank one symmetric spaces, can be viewed as the total
space of a Riemannian submersion, in some cases an orbifold submersion;
see \cite{fz}. The fact that the homogeneous ones also have totally
geodesic fibers motivated A. Weinstein to study Riemannian submersions
with totally geodesic fibers and positive vertizontal curvatures, i.e.,
sectional curvatures of planes spanned by a vector tangent to a fiber and
a vector orthogonal to it (\cite{we}). He called such fiber bundles
{\it fat}, and showed that this much weaker condition already imposes
strong restrictions.

Fat circle bundles are in one to one correspondence with symplectic
manifolds and hence well understood. Therefore, we will restrict
ourselves to bundles whose fiber dimension is bigger than one, which in
turn implies that the dimension of the base must be divisible by 4.

Let $G\to P\overset{\pi}{\longrightarrow} B^{2m}$ be a $G$--principal
bundle with $G$ a compact connected Lie group endowed with a biinvariant
metric $\la\ ,\ \ra$, and $B$ a compact connected manifold. Given a
principal connection~$\theta$ with curvature form $\Omega$, we say that
$\theta$ is {\it $u$--fat} if
$$
\Omega_u= \la\Omega(\cdot\,,\cdot),u\ra \text{ is non-degenerate on } \h,
$$
where $\h $ is the horizontal space of $\pi$ and $u$ lies in the Lie
algebra $\fg$ of $G$. For a connection metric on $P$ (see Section 1) this
condition is equivalent to requiring that the sectional curvatures
spanned by the action field generated by $u$ and any horizontal vector is
positive. In particular, fatness is independent of the metrics on
the base and fiber.

If fatness holds for all $0\neq u\in \fg$ we say that $\theta$ is
{\it fat}, or simply the principal bundle is fat. If it holds for all
$u\neq 0$ in a subset $\fs\subset \fg$, we say that the principal
connection is {\it $\fs$--fat}, or that $\fs$ is fat. Observe that if $u$
is fat, so are all vectors in its adjoint orbit $\fo=Ad_G(u)$. Following
Weinstein, we consider the homogeneous $Ad_G$--invariant polynomial
$q_\fo:\fg\to\R$ defined as
$$
q_\fo(\alpha)=\int_G\la Ad_g(u),\alpha\ra^m dg.
$$
By Chern-Weyl theory, there exists a closed $2m$--form $\omega_\fo$ on
$B^{2m}$ such that $\pi^*\omega_\fo=q_\fo(\Omega)$. By $\fo$--fatness,
the form $\la Ad_g(u),\Omega\ra^m$ is everywhere nonzero on $\h$ and thus
$\omega_\fo$ is a volume form on $B^{2m}$. Therefore, the characteristic
number $\int_B\omega_\fo$ is nonzero and we call it
the {\it Weinstein invariant associated to $\fo$}.

\vskip .2cm

Our main purpose is to compute these invariants for the classical Lie
groups, obtaining explicit topological obstructions to fatness in terms
of Chern (or Pontrjagin) numbers. This will allow us to derive several
applications.
\vskip .2cm

The simplest case is the torus $T^n$ for which we obtain a lower bound on
the Betti numbers of the base.

\begin{main}\label{t:tori}
Let $T^n \to P \to B^{2m}$ be a fat principal bundle. Then, the Betti
numbers of $B^{2m}$ satisfy $b_{2i}\geq n$, for all
$1\leq i\leq m-1$.
\end{main}

Denote by $c_k\in H^{2k}(B,\Z)$ and $p_k\in H^{4k}(B,\Z)$ the Chern and
Pontrjagin classes, and by $e\in H^{2n}(B,\Z)$ the Euler class when
$G=SO(2n)$. In the case of $U(2)$ and $SO(4)$, Weinstein invariants have
rather simple expressions.

\begin{main}\label{t:simpleri}
Let $G\to P\to B^{2m}$ be a fat principal bundle.
\begin{enumerate}[a)]
\item If $G=U(2)$, then $c_1^m$ and $(c_1^2-4c_2)^{m/2}$ are nonzero
and have the same sign, and
$$
\sum_{j=0}^{m/2} \binom{m+1}{2j+1}t^{2j}c_1^{m-2j}(c_1^2-4c_2)^j\neq 0,
\ \ \ \forall \ t\in \R;
$$
\item If $G=SO(4)$, then $(p_1+2e)^{m/2}$ and $(p_1-2e)^{m/2}$ are
nonzero and have the same sign, and
$$
\sum_{j=0}^{m/2} \binom{m+2}{2j+1}t^{2j}(p_1-2e)^{m/2-j}(p_1+2e)^j\neq 0,
\ \ \ \forall \ t\in \R.
$$
\end{enumerate}
\end{main}
\vskip .2cm

For Lie groups with rank$(G)>2$ the formulas become more complicated,
e.g. in the case of $U(n)$ we obtain
$$
\sum_{ n\geq
\lambda_1\geq\dots\geq\lambda_m\geq 0   }\prod_{i=1}^{m}(n+m-i-\lambda_i)!\,
\det\left(\sigma_{\lambda_i+j-i}(y)\right)_{1\leq i,j\leq m}\,
\det\left(c_{\lambda_i+j-i}\right)_{1\leq i,j\leq m}\neq 0
$$
for all $0\ne y=(y_1,\dots , y_n)$, with $\sum\lambda_i=m$, where
$\sigma_{i}$ stands for the elementary symmetric polynomial of degree
$i$ in $n$ variables. Observe that special consequences are $c_1^m\neq 0$
for $y=(1,\cdots,1)$ and
$\det\left(c_{j-i+1}\right)_{1\leq i,j\leq m}\neq0$ for
$y=(1,0,\cdots,0)$. Similar formulas hold for the other classical Lie
groups; see \tref{t:main}.

We use these invariants as follows. Given a $G$--principal bundle $P$,
the Weinstein invariants define homogeneous polynomials in $\rk(G)$
variables once we parametrize the adjoint orbits $\fo\subset\fg$ in terms
of a maximal abelian subalgebra $\ft\subset \fg$, by writing
$\fo=Ad_G(y)$ for $y\in\ft$. The coefficients of these multivariable
polynomials are Chern numbers of $P$, and fatness implies that they have
no nonzero real roots. It is thus in general difficult to express the
nonvanishing of the Weinstein invariants in terms of the Chern numbers
alone. But in some cases this is possible. For example, we have the
following which applies, in particular, to base manifolds with
$b_4(B^{2m})=1$.

\begin{cor}\label{c:cor2i}
Let $G\to P \to B^{2m}$ be a principal bundle, where $G=U(2)$ or
$G=SO(4)$. Suppose there exists $r\in\R$ such that
$c_1^2-4c_2=rc_1^2\neq0$ if $G=U(2)$, or $p_1+2e=r(p_1-2e)\neq0$ for
$G=SO(4)$. Then all Weinstein invariants are nonzero if and only if
$r>0$.
\end{cor}

A natural context where partial fatness appears is for associated
bundles. Given $H\subset G$ a closed subgroup, we have the associated
bundle by $G/H \to P\times_G G/H=P/H \to B$. A connection metric on this
fiber bundle can be described in terms of a principal connection on $P$.
The vertizontal curvatures of such a connection metric are positive if
and only if the principal connection on $P$ is $\fh^\perp$--fat. It turns
out that any fat bundle is associated to some principal bundle in this
way.

Perhaps the most natural examples of associated bundles are the sphere
bundles. In the real case, we conclude from \tref{t:main} the following.

\begin{cor}\label{c:realspheres}
A sphere bundle with totally geodesic fibers and positive vertizontal
curvatures satisfies $\det(p_{j-i+1})_{1\leq i,j\leq\frac{m}{2}}\neq 0$,
where $2m$ is the dimension of the base.
\end{cor}

We will see that if the sphere bundle is the sphere bundle of a complex
or quaternionic vector bundle, we obtain a one parameter family of
obstructions instead of a single one.

\smallskip

In \cite{dr} it was shown that the only $\Sp^3$ bundle over $\Sp^4$ which
admits a fat connection metric is the Hopf bundle. For $\Sp^7$ bundles
over $\Sp^8$ this is still an open problem. Such bundles $\Sp^7\to
M_{k,l}\to \Sp^8$ are classified by two arbitrary integers $k,l$ such
that $p_2=6(k-l)$ and $e=k+l$. Using the obstructions for quaternionic
sphere bundles and \cref{c:realspheres}, we will show
\vskip .4cm

\begin{cor}\label{c:s7}
The sphere bundles $\Sp^7\to M_{k,l}\to \Sp^8$, where either $k=l$ or
$(k,l)=(8r,4r)$, $r\in\Z$, do not admit a fat connection metric.
In particular, for $k=l=1$, it follows that $T_1\Sp^8\to\Sp^8$ does not
admit a fat connection metric.
\end{cor}

Similarly, we will see that $T_1\CP^4\to \CP^4$ admits no fat
connection metric. For fat $\Sp^3$--fiber bundles over $\CP^2$ we have
the following.

\begin{cor}\label{c:zi1}
The only two $\Sp^3$--fiber bundles over $\CP^2$ that may admit a fat
connection metric are the complex sphere bundles with $c_1^2=9$ and
$c_2=1$ or $2$. In particular, $T_1\CP^2\to\CP^2$ does not have a fat
connection metric.
\end{cor}

\vskip .2cm

Bérard--Bergery classified in \cite{bb} the fat fiber bundles which are
homogeneous. A family of such examples are the fiber bundles over the
Grassmannian of 2--planes in $\C^{n+1}$,
$$
U(2)/S^1_{p,q}\to U(n+1)/U(n-1)\cdot S^1_{p,q} \to G_2(\C^{n+1}),
$$
where $ S^1_{p,q}=\{\diag(z^p,z^q)\in U(2) : z \in S^1\}$. The fiber
is the lens space $\Sp^3/\Z_{p+q}$ when $p+q\neq 0$. He showed that
the bundle has a homogeneous fat connection metric if and only if $pq>0$.
We will show that the homogeneity property can be dropped:

\begin{cor}\label{c:lens2} The above
fiber bundles with $pq\leq 0$ admit no fat connection metric.
\end{cor}

For general  $U(2)/S^1_{p,q}$ fiber bundles we have

\begin{cor}\label{c:lens1} If $U(2)/S^1_{p,q}\to P \to B^{2m}$
is a fat bundle, then $(c_1^2-4c_2)^{m/2}\neq 0$. Moreover, if
$(p,q)\neq(1,1)$ and $c_1^2=r(c_1^2-4c_2)$ for some $r\in\R$, then
the nonvanishing of the Weinstein invariants is equivalent to
$r>-\left(\frac{1-\cos(\frac{\pi}{m+1})}{1+\cos(\frac{\pi}{m+1})}\right)
\left(\frac{p+q}{p-q}\right)^2$.
\end{cor}

The above fat bundles over $G_2(\C^{n+1})$ show that \cref{c:lens1} is
sharp, since for $n=2$ these bundles are associated to the same
principal bundle $U(2)\to SU(3)\to \CP^2$ which has $r=-1/3$.

\medskip

In Section 1 we collect various facts about fat fiber bundles, and
the algebra of symmetric polynomials and Schur functions. These turn out
to be central in converting the integrals into polynomials in Chern
and Pontrjagin numbers. In Section 2 we derive the general form of
Weinstein invariants for all classical Lie groups and for $G_2$. In
Section 3 we concentrate on the case of $G=T^n,\ U(2)$ and $SO(4)$, while
in Section 4 we discuss sphere bundles. In Section 5 we prove a stronger
version of the reduction conjecture stated in \cite{zi1} for normal
subgroups, namely, no fat vector exists in $\fh^\perp$ if the structure
group reduces to a normal subgroup with Lie algebra $\fh$.
Finally, in Section 6 we relate our obstructions to some of the fat
bundles in \cite{bb} by computing their Weinstein invariants. Throughout
the paper we will provide several additional applications.

\vskip .2cm
We would like to thank N. Wallach for helpful conversations.

\section{Preliminaries}  

We first recall Weinstein's definition of fatness of a fiber bundle and
his basic topological obstruction to fatness; see \cite{we} and
\cite{zi1}.

\bigskip

Let $\pi\colon M\to B$ be a fiber bundle with fiber $F$, and metrics on
$M$ and $B$ such that $\pi$ is a Riemannian submersion. Let $\h$
and $\v$ denote the horizontal and vertical subbundles of $TM$. If the
fibers of $\pi$ are totally geodesic, the sectional curvature of a
vertizontal 2-plane, i.e. a plane spanned by a vertical vector $U$ and
horizontal vector $X$ is equal to $\|A_UX\|^2$, where $A\colon
\h\times\h\to\v$ is the O'Neill tensor and $\ml A_UY, X\mr=-\ml
A_XY,U\mr$. In particular, these curvatures are automatically
non--negative. According to Weinstein,
$$
\pi\colon M \to B \ \ {\rm is\ called\ {\it fat}\ if }\
A_XU\ne 0 \ \, {\rm for\ all }\ 0\ne X\in\h,\ 0\ne U\in\v,
$$
or, equivalently, when all vertizontal sectional curvatures are positive.

\medskip

We first consider the case where the fiber bundle $\tau\colon P\to B$ is
a $G$--principal bundle and the horizontal spaces are $G$--invariant. The
horizontal distribution can then be described in terms of a principal
connection $\theta\colon TP\to \fg$ as $\h=\ker\theta$, where $\fg$ is
the Lie algebra of $G$. With the aid of a metric on the base and a left
invariant metric on~$G$, $\theta$ defines a so called {\it connection
metric} on $P$ by declaring $\h$ and $\v$ to be orthogonal, endowing $\h$
with the pull back of the metric on the base, and $\v$ with the chosen
left invariant metric on $G$. If we endow $P$ with such a connection
metric and $\Omega\,\colon\,TP\times TP\to\fg$ is the curvature form
of $\theta$, Weinstein observed that fatness of $\tau$ can be
rewritten as:
\begin{equation}\label{e:fatpb}
{\rm For\ each}\ \, 0\neq y\in\fg,\
\Omega_y:=\la\Omega(\cdot\,,\cdot),y\ra\ {\rm is\ a\ non\!-\!degenerate\
\, 2\!-\!form\ on}\ \h,
\end{equation}
where we have chosen an auxiliary biinvariant metric $\la\cdot,\cdot\ra$
on~$G$. This is indeed an immediate consequence of
$2\theta(A_XY)=-\Omega(X,Y)$. In particular \eqref{e:fatpb} implies that
fatness is independent of the metrics on the base and fiber, i.e., it
only depends on the principal connection itself. We thus simply say that
the principal connection $\theta$, or $P$ by abuse of language, is
{\it fat}. Furthermore, if $\fs\subseteq\fg$ is a subset, we will say
that $\theta$ is {\it $\fs$--fat} if $\Omega_y$ is non--degenerate for
all $0\neq y\in\fs$. Also observe that if a vector $y\in\fg$ is fat, the
whole adjoint orbit $Ad_G(y)$ consists of fat vectors since
$\Omega_{Ad_g(y)}=g^*(\Omega_y)$. Hence we can assume
that $y$ lies in a maximal abelian subalgebra $\ft\subseteq\fg$.

Observe that for each fat vector $y\in\fg$ we have a nonvanishing vector
field $Z_y$ on the unit sphere of $\h$ given by
$\la Z_y(X),Y\ra= \Omega_y(X,Y)$, and if $\{y_1,\dots,y_r\}$ is a basis
of a fat subspace $V\subset \fg$, the vector fields
$Z_{y_1},\dots,Z_{y_r}$ are pointwise linearly independent. By the well
known Radon--Hurwitz formula $V$--fatness thus implies:
\begin{equation}\label{e:dim}
\text{If } \dim B=(2a\!+\!1)2^{4b+c} \
\text{ with $0\le c\leq 3$, then }
\dim V\leq 2^c+8b-1.
\end{equation}
In particular,
\begin{equation}\label{e:fatdim}
{\rm If}\ \dim V \ge 2, 4, 8,\ {\rm then}\ \ 4|\dim B,\ \ 8|\dim B,\ \
16|\dim B,\ {\rm respectively.}
\end{equation}
Notice that the adjoint orbit of $V$ may contain a linear subspace of
larger dimension, in some cases all of $\fg$, which gives further
restrictions.

\medskip

We define the Weinstein invariants as follows, where we assume
that $G$ and $B$ are compact and connected. For each adjoint
orbit $\fo\subseteq\fg$, we write $\fo=Ad_G(y)$ for $y\in\ft$. For
$k\in\N\cup\{0\}$, the homogeneous $Ad_G$--invariant polynomial
$q^k_y=q^k_\fo:\fg\to\R$ given by
\begin{equation}\label{e:main}
q^k_\fo(\alpha)=\int_G\la Ad_g(y),\alpha\ra^k dg
\end{equation}
defines a closed $2k$--form $\omega_\fo$ on $B^{2m}$ via
$\tau^*\omega_\fo=q^k_\fo(\Omega)$. By Chern--Weyl theory,
$[w_\fo]\in H^{2k}(B,\R)$ represents a characteristic class of the
bundle. Now
suppose that $k=m$ is half the dimension of the base and write
$q_\fo=q^m_\fo$. If $\fo$ is fat, $\Omega_y^m\ne 0$ is a volume form on
$\h$. Thus, if $G$ is connected,
$\la Ad_g(y),\Omega\ra^m$ is nowhere zero and has constant sign when $g$
varies along $G$, and the integral $q_\fo(\Omega)$ is nonzero on $\h$.
Hence $\omega_\fo$ is a volume form of $B^{2m}$, in particular $B^{2m}$
is orientable, and the characteristic number $\int_B \omega_\fo$ is
nonzero. We call this characteristic number the {\it Weinstein invariant
associated to $\fo$}, and our main goal is to express it explicitly in
terms of Chern and Pontrjagin numbers. This will allow us to obtain
various applications.

Observe that, for a circle bundle, fatness is equivalent to $\omega$
being a symplectic form on the base, where $\omega$ is given by
$\Omega=\tau^*w$. Thus, the only Weinstein invariant is the symplectic
volume. For any other fat fiber bundle, by \eqref{e:fatdim} we have that
4 divides the dimension of the base. Therefore, we always assume
that $m$ is even.

\medskip

For convenience, we use the same notation for an $Ad_G$--invariant
polynomial on $\fg$, for its restriction to a maximal abelian subalgebra
$\ft\subseteq\fg$, and for the corresponding characteristic class.
Recall that the {\it Chern classes} $c_i\in H^{2i}(B,\Z)$ are defined by
the $Ad_{U(n)}$--invariant polynomials $c_i(A)$, with
$$
\det(I+tA)=\sum_ic_i(A)t^i\ ,\ \ A\in\fu(n),
$$
and the {\it Pontrjagin classes} $p_i\in H^{4i}(B,\Z)$ by the
$Ad_{O(n)}$--invariant polynomials $p_i(A)$, with
\begin{equation}\label{e:pont}
\det(I+tA)=\sum_ip_i(A)t^{2i}\ ,\ \ A\in\fo(n).
\end{equation}
For even rank, we also have the {\it Euler class} $e\in H^{2n}(B,\Z)$
given by the $Ad_{SO(2n)}$--invariant Pfaffian
$$
e(A)=\Pf(A),\ \  A\in\fo(2n),
$$
where $\Pf^2(A)=\det(A)$. Again, by abuse of notation, we also use the
same symbol for the {\it quaternionic Pontrjagin classes}
$p_i\in H^{4i}(B,\Z)$ given by the $Ad_{Sp(n)}$--invariant polynomials
$p_i(A)$ as in \eqref{e:pont}, but for $A\in\fsp(n)$.
As it is well known, these polynomials form a basis of the set of all
$Ad_G$--invariant polynomials in the case of a classical Lie group
$G$, the only relations being $e^2=p_n$ in the case of $G=SO(2n)$ and
$c_1=0$ for $G=SU(n)$. Thus each Weinstein invariant is a polynomial in
these basic classes, $c_i,\ p_i$ and $e$, evaluated on the fundamental
cycle $[B]$.

Now, as a function of $\fo=Ad_G(y)$, $q_\fo=q_y$ becomes a polynomial in
$y\in\ft\cong\R^n$, with coefficients being Chern or Pontrjagin numbers.
By definition they are invariant under the Weyl group $W=N(T)/T$. We use
the following standard forms for $\ft$:
for $\fu(n)$ and $\fsp(n)$ we have
$\ft=\{i\,\diag(y_1,\dots,y_n):y\in\R^n\}$, while for $\fs\fo(2n+1)$ and
$\fs\fo(2n)$ we have $\ft=\{\diag(y_1J,\dots,y_nJ):y\in\R^n\}$, where $J$
stands for the basic $2\times 2$ skew symmetric matrix. We will denote
both the vector and its coordinates by $y=(y_1,\cdots,y_n)\in \ft$. Since
all Weyl groups of the classical Lie groups contain the permutation
group, $q_y$ can be expressed in terms of the
{\it elementary symmetric polynomials}
$\sigma_i(y)=\sigma_i(y_1,\dots,y_n)$.
We will choose the biinvariant metric on $G$ in such a way that the
canonical basis in $\R^n\cong \ft$ is orthonormal.

\bigskip

To obtain the invariants it was important for $G$ to be a connected
compact Lie group. We claim that the obstructions also hold for
non--connected groups. To see this we lift the bundle to a certain cover
of the base whose structure group reduces to a connected Lie group. Let
$\Gamma=G/G_o$ be the component group of the Lie group, on which $G$ acts
naturally. Let $\bar{B}$ be a connected component of
$\{(b,\gamma)\mid b\in B,\ \gamma\in\Gamma\}$ and
define the cover $\alpha \colon \bar{B} \to B$ by $\alpha((b,\gamma))=b$.
This induces the pull back bundle
$\alpha^*(P)=\{(x,b,\gamma)\in P\times\bar{B}\mid \sigma(x)=b\}$
on which $G$ acts via $g\*(x,b,\gamma)=(xg^{-1},b,g\gamma)$. There now
exists a reduction $\bar{P}=\{(x,b,e)\}\subset \alpha^*(P)$ which is
preserved by $G_o$. The connection $\theta$ on $P$ pulls back to a
connection $\bar{\theta}$ on $\bar{P}$ and its curvature $\bar{\Omega}$
still satisfies the property that $\bar{\Omega}_y$ is a non-degenerate
2-form on $\h$ for all $0\neq y\in\fg$. Thus, if
$\bar{q}_\fo(\alpha)=\int_{G_o}\la Ad_g(y),\alpha\ra^m dg$, we have
$\bar{q}_\fo(\bar{\Omega})=\bar{\tau}^*(\bar{\omega}_\fo)$
with $\bar{\omega}_\fo$ a volume form on $\bar{B}$, and
$\bar\omega_\fo=\alpha^*(\omega_\fo)$. Therefore, $\omega_\fo$ is also
a volume form and $\int_{{B}}{\omega}_\fo\ne 0$.

\medskip

We now study how Weinstein invariants behave for coverings.
Let $\tilde G$ be a
finite cover of a connected Lie group $G$, $G=\tilde G/\Gamma$, and
assume that $P$ has a cover $\varphi:\tilde P\to P$ which is a $\tilde
G$--principal bundle. Due to the fact that $\Gamma$ is a subgroup of
the center of $G$, \eqref{e:main} is invariant under $\Gamma$ and
therefore the Weinstein invariants for $\tilde P$ are precisely the ones
for $P$ multiplied by the order of $\Gamma$. Moreover, observe that a
principal connection $\theta$ on $P$ is fat if and only if
$\varphi^*\theta$ on $\tilde P$ is fat.

Similarly, if $G\to G/\Gamma=G^*$ is a covering, a $G$ principal bundle
$P\to B$ induces a $G^*$ principal bundle $P/\Gamma\to B$. If $\theta$ is
a fat connection on $P$, there exists a connection $\theta^*$ on
$P/\Gamma$ whose pullback is $\theta$ since $\Gamma$ is a subgroup of the
center. Again, $\theta^*$ is fat if $\theta$ is fat and the Weinstein
invariants are the same up to a constant.

\medskip

Finally, if $G$ is a product group, $G=G_1\times G_2$, or a local product
$G=G_1\cdot G_2=(G_1\times G_2)/\Gamma$, then from \eqref{e:main} we
get, up to a factor,
\begin{equation}\label{e:peoduct}
q_{(y_1,y_2)}(\alpha_1,\alpha_2)=
\sum_{i=0}^{m}\binom{m}{i} q^i_{y_1}(\alpha_1)q^{m-i}_{y_2}(\alpha_2),
\end{equation}
where $q^k_{y_j}$ are given by \eqref{e:main} for each $G_j$ (we know
they are nonzero only for $k=m$). In particular, if we change the
biinvariant metric on $G$ by multiplying by a constant $c_i$ on each
factor $G_i$, the Weinstein invariants change by a constant as well,
once we replace $(\alpha_1,\alpha_2)$ by $(c_1\alpha_1,c_2 \alpha_2)$,
and its nonvanishing is thus independent of the choice of biinvariant
metrics.
\bigskip

We now discuss the case of a fiber bundle $\pi\colon M\to B$ with fiber
$F$, where we allow a general Riemannian submersion with totally geodesic
fibers. The fiber bundle $\pi$ is associated to a $G$--principal bundle
$\tau\colon P\to B$\, via $M=P\times_GF$, where $G$ acts on $P$ on the
right and on $F$ on the left, $[(p,h)]=\{(pg,g^{-1}h):g\in G\}$, and
$\pi$ can be regarded as the projection onto the first factor. Choose a
principal connection $\theta:TP\to\fg$, a metric on $F$ invariant under
the action of $G$ and a metric on $B$. The horizontal space at $p\in P$
given by $\h_p=\ker\theta_p$ defines a horizontal space at $x=[(p,h)]\in
M$ via $\h_{[(p,h)]}=[(\h_p,0)]$. We now define a metric on $M$ by
pulling back the metric on $B$ with $\pi$, by declaring the fibers to be
orthogonal to $\h$, and choosing the metric on $\pi^{-1}(b)\simeq F\colon
p\to [(p,h)],\ p\in \tau^{-1}(b)$, to be the given metric on~$F$. In this
metric, $\pi$ is a Riemannian submersion with totally
geodesic fibers isometric to $F$, and any Riemannian submersion with
totally geodesic fibers arises in this fashion for some principal bundle;
see \cite{zi1} and \cite{gw} for details. Notice also that, in contrast
to \cite{we}, this metric does not require a choice of metrics on $G$ or
$P$. The metrics described as above are often called {\it connection
metrics} of the fiber bundle.

\medskip

Bérard--Bergery showed in \cite{be} that the holonomy group of a fat
fiber bundle acts isometrically and transitively on the fibers. Since the
holonomy group is contained in $G$, $G$ acts transitively on the fibers
as well. Hence we can assume that $F={G}/{H}$ for some subgroup ${H}$ and
$M=P\times_GG/H=P/H$.
If $\fh^\perp\subset\fg$ is the orthogonal complement of the Lie algebra
$\fh$ of $H$ with respect to our fixed auxiliary
biinvariant metric on $G$, Weinstein showed that:
$$
{\rm The\ connection\ metric\ on}\ \pi\colon M=P\times_GG/H\to B\
{\rm is\ fat\ if\ and\ only\ if}\ \Omega\ {\rm is}\ \fh^\perp{\rm -fat},
$$
i.e. $\Omega_u$ is non-degenerate on $\h$ for all $0\neq u\in
\fh^\perp$. Notice that this is again independent of the $G$--invariant
metric on $F=G/H$ and the metric on $B$, in other words, it only depends
of the principal connection. Therefore, fatness of $\pi$ implies the
nonvanishing of the Weinstein invariants of the $G$--principal bundle
associated to any $0\neq y\in \fh^\perp$, and we write these in
terms of the characteristic numbers of the $G$--principal bundle.

\smallskip

If the left action of $G$ on $G/H$ extends to an action of $G'$ with
$G/H=G'/H'$, one can view, as above, the metric on $M$ induced by
$\theta$ as the metric induced by the unique extension of the principal
connection $\theta$ to a connection $\theta'$ on $P'=P\times_GG'$.
Furthermore, it follows that if $\theta$ is fat then $\theta'$ is fat as
well. Indeed, on $P\subset P'$, $\theta= \theta'|_P$ and thus
$\Omega=\Omega'|_P$, in particular, $\Omega'|_P$ has values in
$\fg\subset\fg'$ and using the identification $TF=\fg^\perp\simeq
\fg'^\perp$ it follows that $\Omega$ is $\fh^\perp$--fat if and only if
$\Omega'|_P$ is $\fh'^\perp$--fat. Thus the possible principal bundles
are all extensions of the holonomy bundle and the principal connection is
the uniquely defined extension of the tautological principal connection
of the holonomy bundle. Nevertheless, the Weinstein invariants depend on
the particular choice of the principal bundle, a fact that we will be
able to exploit in certain situations.

\vskip .2cm

\subsection{Symmetric polynomials and Schur functions}

Fix a positive integer $n\in\N$. A base of the algebra of symmetric
polynomials in $x=(x_1,\dots,x_n)$ is indexed by {\it partitions}
$\lambda=(\lambda_1,\dots,\lambda_n)\in\N_0^n=(\N\cup\{0\})^n$, where
$\lambda$ is non-increasing, i.e., $\lambda_i\geq\lambda_{i+1}$. Denote
by
$$
K_m=\{\lambda\in\N_0^n:\lambda_1\geq\dots\geq\lambda_n,\ |\lambda|=m\}
$$
the set of partitions of degree $m$, where for each $\mu\in\N_0^n$ we set
$|\mu|=\sum_{j=1}^n \mu_j$. We also set $\mu!=\mu_1!\cdots\mu_n!$ and
$x^\mu=x_1^{\mu_1}\cdots x_n^{\mu_n}$ for $x\in\R^n$. We use the
convention $K_m=\emptyset$ if $m$ is not an integer.
For $k\in\N_0$ we also denote by $k=(k,\dots,k)\in\N_0^n$ and
$k\mu=(k\mu_1,\dots,k\mu_n)$. The notation $\mu\subset \gamma$ means that
$\mu_i\leq \gamma_i$ for all $i$. We say that $\mu$ is {\it even} (resp.
{\it odd}) if each $\mu_i$ is even (resp. odd). The partition $\lambda'$
{\it conjugate} to $\lambda\in K_m$ is the partition $\lambda'\in\N_0^m$
defined as $\lambda'_i=\#\{j:\lambda_j\geq i\}$, $1\leq i\leq m$. Since
$\lambda'\subset n$, the set of conjugate partitions to $K_m$ is
$$
K'_m=\{\lambda\in\N_0^m:n\geq \lambda_1\geq\dots\geq\lambda_m,
\ |\lambda|=m\}.
$$
A basic property of the conjugated partitions is that
$\lambda''=\lambda$ for $\lambda \in K_m$ and hence $K_m^{''}=K_m$.
For $\lambda \in K'_m$, we denote by $n-\lambda$
the partition $(n-\lambda_m,n-\lambda_{m-1},\dots,n-\lambda_1)$.

\medskip

Associated with each $\mu\in\N_0^n$ there
is an {\it alternant} $A_\mu$ defined by
\begin{equation}\label{e:das}
A_\mu(x)=\det(x_i^{\mu_j})=\sum_{\sigma\in S_n}\sg(\sigma)x^{\sigma \mu}.
\end{equation}
By definition one has $A_{\sigma\mu}=\sg(\sigma)A_\mu$ for all $\sigma$
in the permutation group $S_n$ of $n$ elements. Since $A_\mu=0$ when
$\mu$ has repeated indexes, a nonvanishing $A_\mu$ can be written, up to
sign, as $A_{\lambda+\rho}$ for $\lambda\in K_m$, where
$m=|\mu|-n(n-1)/2$ and
$$
\rho=\rho_n:=(n-1,n-2,\dots,1,0).
$$
As a special case we have the Vandermonde determinant
$$
\Delta(x):=A_\rho(x)=\det(x_i^{j-1})=\Pi_{i<j}(x_i-x_j).
$$
Since, for any partition $\lambda$, $A_{\lambda+\rho}(x)=0$ if
$x_i=x_j$ for some $i\neq j$, we have that $\Delta$ divides
$A_{\lambda+\rho}$. This allows us to define the degree $m$
homogeneous symmetric {\it Schur polynomial}
$$
S_\lambda:=A_{\lambda+\rho}\,/\Delta,\ \ \ \ \ \ \lambda\in K_m.
$$

\medskip

The {\it complete homogeneous symmetric polynomial} of degree $m$ in $n$
variables is defined as
$$
h_m=h_m(x) := \sum_{|\mu|=m} x^\mu,
$$
and we set $h_\lambda=\Pi_ih_{\lambda_i}$ for any partition $\lambda$.
The {\it elementary symmetric polynomials} are
$$
\sigma_m:=\sum_{|\mu|=m,\ \mu\subset 1} x^\mu.
$$
There is a simple relation between these polynomials given by
\begin{equation}\label{e:inv}
h_0=1,\ \ \ \ \ \sum_{j=0}^r (-1)^j\sigma_jh_{r-j}=0,
\ \ \ \forall \ r\geq 1,
\end{equation}
where we set $\sigma_m=0$ for $m>n$ or $m<0$ (cf. (2.6$'$) in
\cite{ma}, or \cite{fh}). It follows that $\{h_0,h_1,\dots,h_n\}$ is
another basis of the algebra of symmetric polynomials. Moreover, one has
$$
h_m=\det(\sigma_{j-i+1})_{1\leq i,j\leq m},\ \  \forall m\in\N
$$
(cf. page 20 in \cite{ma}). For the special cases $n=2$ or
$x_3=\dots=x_n=0$ we also have that
\begin{equation}\label{e:h2}
h_m=\frac{1}{2^m}\sum_{j=0}^{[m/2]}\binom{m+1}{2j+1}
\sigma_1^{m-2j}(\sigma_1^2-4\sigma_2)^j.
\end{equation}
Indeed, by definition, $h_m=(x_1^{m+1}-x_2^{m+1})/(x_1-x_2)$.
If we consider $x_1$ and $x_2$ formally as the roots of
$x^2-\sigma_1 x +\sigma_2=0$ and apply the binomial theorem
one obtains \eqref{e:h2}.

\begin{remark}\label{r:neg}
{\rm Because of \eqref{e:inv}, when $E$ is a complex vector bundle, $h_m$
represents the $m$--Chern class $c_m(-E^*)$ of the formal negative of the
dual of $E$.}
\end{remark}

We summarize now the main properties of Schur functions that will be
needed.

\begin{proposition}\label{p:si}
We have the following properties of Schur functions
$S_\lambda(x_1,\dots,x_n)$ for $\lambda \in K_m$:
\begin{enumerate}[{\rm (}a\,{\rm )}]
\item $S_{\lambda}=\det(h_{\lambda_i+j-i})_{1\leq i,j\leq n}$;
\item $S_{\lambda}=\det(\sigma_{\lambda'_i+j-i})_{1\leq i,j\leq m}$,
where $\lambda'\in K'_m$ is the conjugate partition of $\lambda$;
\item $S_{\lambda}=\sigma_m$, for $\lambda=(\underbrace{1,\dots,1}_m,0,\dots,0)$;
\item $S_{\lambda}=h_m=\det(\sigma_{j-i+1})_{1\leq i,j\leq m}$,
for $\lambda=(m,0,\dots,0)$;
\item $S_\lambda=\sigma_n^{\lambda_n}S_{\tilde\lambda}$,
where $\tilde\lambda=\lambda-\lambda_n1$;
\item $S_0=1$;
\item $S_\lambda(e_1)=1$, for $\lambda=(m,0,\dots,0)$ and 0 for any other
partition $\lambda$;
\item $S_\lambda(e_1+te_2)=\sum_{i=k}^{m-k}t^i$,
if $\lambda=(m\!-\!k,k,0,\dots,0)$, $0\leq k \leq [m/2]$,
and 0 for any other partition $\lambda$.
\end{enumerate}
\end{proposition}
\proof
All statements are immediate consequences of the well known Jacobi-Trudi
identities $(a)$ and $(b)$ (see (3.4) and (3.5) in \cite{ma}).
\qed
\vspace{1.5ex}

We will see in the next section that Weinstein invariants can be
explicitly written in terms of Schur functions and then, by part $(b)$
of \pref{p:si}, in terms of Chern and Pontrjagin numbers.

\section{Weinstein invariants for the classical groups} 

Assume $G$ is connected. To convert the integral in \eqref{e:main} into
a sum, we use a formula due to Harish--Chandra (see Theorem 3.2.1.3 in
\cite{wa})
$$
t^r\pi(y)\pi(x)\int_Ge^{t\la Ad_g(y),x\ra}=
\sum_{w\in W}\det(w)\,e^{t\la wx,y\ra}, \ \ \ \ t\in\R,
$$
which holds after multiplying the biinvariant metric with a suitable
constant. Here, $W$ denotes the Weyl group of $G$, $\pi$ the product of
the positive roots of $(\fg,\ft)$, $r$ the number of positive roots,
i.e., $2r=\dim G - \rk (G)$, and
$x,y\in\ft\cong\R^n$. Comparing Taylor series expansions in $t$, we have
from \eqref{e:main} that
\begin{eqnarray*}
\pi(x)\pi(y)q^k_\fo(x)&=&
\frac{k!}{(k+r)!}\sum_{w\in W}\det(w)\la wx,y\ra^{k+r}\\
&=&k!\sum_{w\in W}\det(w)\sum_{|\mu|=k+r}\frac{1}{\mu!}(wx)^\mu y^\mu\\
&=&\sum_{|\mu|=k+r}\frac{k!}{\mu!}\ y^\mu\sum_{w\in W}\det(w)(wx)^\mu,
\end{eqnarray*}
with $k\in\N$, $\mu\in\N_0^n$. Therefore,
\begin{equation}\label{e:invl}
\pi(x)\pi(y)q^k_\fo(x)=\sum_{|\mu|={k+r}}\frac{k!}{\mu!}\ y^\mu\ L_\mu(x),
\end{equation}
where
$$
L_\mu(x):=\sum_{w\in W}\det(w)(wx)^\mu.
$$

\bigskip

We now compute the polynomials $q^k_y$, $y\in\ft$, for the
classical groups in terms of Schur functions. Since they are
$Ad_G$--invariant we only need to describe their restriction to $\ft$.
Recall that $K_r$ is empty if $r$ is not an integer, and we set a sum
over the empty set to be zero.

\begin{proposition}\label{p:main}
Let $G$ be one of the classical groups with Lie algebra $\fg$ and
$\ft\cong\R^n$ the Lie algebra of a maximal torus of $G$. Then, up to a
positive constant which only depends on $G$, we have the following
expressions for $q^k_y:\ft\to\R$ for each $y\in\ft$ and $k\in\N\!:$
\begin{enumerate}[a)]
\item If $G=U(n)$ or $SU(n)$, then
$$
q^k_y(x)=\sum_{\lambda\in K_k}\!\!
\frac{k!}{(\lambda\!+\!\rho)!}S_\lambda(y)S_\lambda(x);
$$
\item If $G=O(2n),O(2n+1),SO(2n+1)$ or $Sp(n)$, then
$$
q^k_y(x)=\sum_{\lambda\in K_{k/2}}
\frac{k!}{(2(\lambda\!+\!\rho)\!+\!\epsilon)!}
S_\lambda(y^2)S_\lambda(x^2),
$$
where $\epsilon=0$ for $G=O(2n)$ and $\epsilon=1$ otherwise;
\item If $G=SO(2n)$, then
$$
q^k_y(x)=\!\!\!\!\sum_{\lambda\in K_{k/2}}\!\!\!
\frac{k!}{(2(\lambda\!+\!\rho))!}S_\lambda(y^2)S_\lambda(x^2)+
\!\!\!\!\!\!\!
\sum_{\lambda\in K_{(k-n)/2}}\!\!
\frac{k!}{(2(\lambda\!+\!\rho)\!+\!1)!}
\sigma_n(y)S_\lambda(y^2)\sigma_n(x)S_\lambda(x^2).
$$
\end{enumerate}
\end{proposition}
\proof
We proceed case by case since the actual expressions in
\eqref{e:invl} involve the structure of the Lie algebra of $G$.

\vspace{1.5ex}

$G=U(n)$. The Weyl group $W=S_n$ acts on $\ft$ as the permutation
group $S_n$ of $n$ elements. Furthermore, $r=n(n-1)/2$,
$\pi(x)=\Delta(x)$ and $L_\mu(x) = A_\mu(x)$. Since $L_\mu=0$ if there
are two repeated integers in $\mu$ we obtain from
\eqref{e:invl} that
$$
\pi(x)\pi(y)q^k_y(x)=\!\sum_{\mu-\rho\in K_k}\!\frac{k!}{\mu!}
\sum_{\tau\in S_n}\ y^{\tau \mu}L_{\tau \mu}(x)=
\!\sum_{\mu-\rho\in K_k}\!\frac{k!}{\mu!}A_\mu(y)A_\mu(x),
$$
which proves \pref{p:main} $(a)$ for the unitary group.

\vspace{.5ex}

$G=SU(n)$. Identifying the maximal torus of $SU(n)$ with
$\{x\in\R^n:\sigma_1(x)=0\}$, the same formula as for $U(n)$ holds
simply taking into account that $\sigma_1(x)=\sigma_1(y)=0$, since the
Weyl group and the roots of $SU(n)$ coincide with those of $U(n)$.

\vspace{1.5ex}

$G=SO(2n)$. Here, $W=S_n\times\Z_2^{n-1}$ acts on $\ft$ as the
permutation group and by an even change of signs, and $r=n(n-1)$. Observe
that $L_\mu=0$ if $\mu$ contains both an even and an odd index or when
two indexes are repeated, and $L_\mu=2^{n-1}A_\mu$ otherwise.
On the other hand we have
$L_{2\mu}(x)=L_\mu(x^2)=2^{n-1}A_\mu(x^2)$ and
$L_{2\mu+1}(x)=\sigma_n(x)L_\mu(x^2)=2^{n-1}\sigma_n(x)A_\mu(x^2)$.
Therefore,
\begin{eqnarray*}
\pi(x)\pi(y)q^k_y(x)&=&\sum_{|\mu|=k+r,\, \mu {\rm \, even\, or\, odd}}
\frac{k!}{\mu!}\sum_{\tau\in S_n}\ y^{\tau \mu}L_{\tau \mu}(x)\\
&=&\!\sum_{\mu-\rho\in K_{k/2}}\!
\frac{k!\,2^{n-1}}{(2\mu)!}A_\mu(y^2)A_\mu(x^2)\\
& &\ \ +\sigma_n(x)\sigma_n(y)\!\!\!\sum_{\mu-\rho\in K_{(k-n)/2}}\!
\frac{k!\,2^{n-1}}{(2\mu\!+\!1)!}A_\mu(y^2)A_\mu(x^2).
\end{eqnarray*}
We immediately get \pref{p:main} $(c)$ up to a factor $2^{n-1}$
since for $SO(2n)$ we have $\pi(x)=\Delta(x^2)$.

$G=SO(2n+1)$. Here, $W=S_n\times\Z_2^n$ acts on $\ft$ as the permutation
group and by arbitrary sign changes, and $r=n^2$. Thus, $L_\mu=0$ if
$\mu$ contains an even index, and
$L_{2\mu+1}(x)=\sigma_n(x)L_\mu(x^2)=2^n\sigma_n(x)A_\mu(x^2)$.
Since $\pi(x)=\sigma_n(x)\Delta(x^2)$, we get
\pref{p:main} $(b)$ for $SO(2n+1)$ up to a factor $2^n$.

\vspace{1.5ex}

$G=O(2n+\e)$, $\e=0,1$.
These groups share a maximal torus with $SO(2n+\e)$. However, we cannot
apply Harish--Chandra's formula directly, since the orthogonal group is
not connected. Write $G=G_0\cup g'G_0$ where \mbox{$G_0=SO(2n+\epsilon)$}
and $g'={\rm diag}\ (-1,1,\dots,1)\in G\setminus G_0$. From
\eqref{e:main} we get
$$
q^k_y(x)=\int_G\la Ad_g(y),x\ra^k dg=
\hat q^k_y(x)+\hat q^k_y\left(Ad_{g'}(x)\right),
$$
where $\hat q^k_y(x)=\int_{G_0}\la Ad_g(y),x\ra^k dg$. Since $Ad_{g'}$
preserves $S_\lambda(x^2)$ and changes the sign of $\sigma_n(x)$,
\pref{p:main} $(b)$ and $(c)$ for $SO(2n+\e)$ imply
\pref{p:main} $(b)$ for the orthogonal group up to a factor $2^{n+\e}$.

\vspace{1.5ex}

$G=Sp(n)$. The symplectic group $Sp(n)$ shares with $U(n)$ the same
maximal torus. The Weyl group acts on it in the same way
as the one of $SO(2n+1)$, but $\pi(x)=2^n\sigma_n(x)\Delta(x^2)$.
Therefore, the expression differs from the one for $SO(2n+1)$ only by a
$2^n$ factor.
\qed

\bigskip
We now have all the ingredients needed to express Weinstein invariants
in terms of Chern and Pontrjagin numbers. Recall that for any fat fiber
bundle where the fibers are different from $\Sp^1$ the dimension of the
base must be divisible by 4.

\begin{theorem}\label{t:main}
Let $G$ be one of the classical groups and $\ft\cong\R^n$ a maximal
abelian subalgebra of its Lie algebra. Let $G\to P \to B^{2m}$ be
a $G$--principal bundle and an element $y\in\ft$ that is fat. Taking into
account that in all statements the indexes $i,j$ of the matrices run
over $1\leq i,j\leq r$ for $\lambda \in K'_r$, we have:
\begin{enumerate}[a)]
\item If $G=T^n$, then
$$
\left(\sum_{i=1}^ny_ic_i\right)^m\neq 0,
$$
where $c_1,\dots,c_n\in H^2(B,\Z)$ are the Chern classes of $P$;
\item If $G=U(n)$ or $SU(n)$, then
$$
\sum_{\lambda\in K'_m}(n\!-\!\lambda\!+\!\rho_m)!\,
\det\left(\sigma_{\lambda_i+j-i}(y)\right)\,
\det\left(c_{\lambda_i+j-i}\right)\neq 0,
$$
where $c_k\in H^{2k}(B,\Z)$ is the $k^{th}$--Chern class of $P$,
with $\sigma_1(y)=0$ and $c_1=0$ for $G=SU(n)$;
\item If $G=O(2n),O(2n+1),SO(2n+1)$ or $Sp(n)$, then
$$
\sum_{\lambda\in K'_{m/2}}(2(n\!-\!\lambda\!+\!\rho_{m/2})\!+\!\epsilon)!\,
\det\left(\sigma_{\lambda_i+j-i}(y^2)\right)\,
\det\left(p_{\lambda_i+j-i}\right)\neq 0,
$$
where
$p_k\in H^{4k}(B,\Z)$ is the $k^{th}$--Pontrjagin class of $P$,
with $\epsilon=0$ for $G=O(2n)$ and $\epsilon=1$ otherwise;
\item If $G=SO(2n)$, then
\begin{eqnarray*}
&{}&\!\!\!\!\!\!\!\!\!\!\!\!\!\!\!\!\!\!\!\!\!\!\!\!
\sum_{\lambda\in K'_{m/2}}
\frac{2(n\!-\!\lambda\!+\!\rho_{m/2}))!}{(2\rho_{(m+2n)/2})!}
\det\left(\sigma_{\lambda_i+j-i}(y^2)\right)
\det\left(p_{\lambda_i+j-i}\right)\\
\ \ \ \ \ \ &+&
e\!\!\!\!\!\sum_{\lambda\in K'_{(m-n)/2}}\!\!\!\!
\frac{(2(n\!-\!\lambda\!+\!\rho_{(m-n)/2})\!+\!1)!}{(2\rho_{(m+n)/2}+1)!}
\sigma_n(y)\det\left(\sigma_{\lambda_i+j-i}(y^2)\right)
\det\left(p_{\lambda_i+j-i}\right)\neq 0,
\end{eqnarray*}
where $e\in H^{2n}(B,\Z)$ is the Euler class of $P$.
\end{enumerate}
\end{theorem}

\proof
For the torus, the Weyl group $W$ is trivial and we
simply get, using Fubini and \eqref{e:main}, that
$$
q_\fo(x)=\left(\sum_{i=1}^ny_ix_i\right)^m,
$$
which gives us part $(a)$.

Using that $\lambda''=\lambda \in K_m$ and \pref{p:si} $(b)$,
the other cases are direct consequences of \pref{p:main} for $k=m$,
writing the expressions in terms of conjugate partitions, and using that
$\rho_{m+n}!=(\lambda+\rho)!(n-\lambda'+\rho_m)!$ and hence
$(2\rho_{m+n}+\e)!=(2(\lambda+\rho)+\e)!\,(2(n-\lambda'+\rho_m)+\e))!$.
Indeed, this follows from the fact that
$$
\{\lambda_i+\!n-i: 1\leq i\leq n\}\cup
\{n+j-\lambda'_j-1: 1\leq j\leq m\}=\{0,1,2,\dots,m+n-1\},
$$
with the union being disjoint; see (1.7) in \cite{ma}.
\qed

\begin{remark}\label{r:g2} Weinstein invariants for $G=G_2$.
{\rm Our methods apply to all Lie groups, and not only to the classical
ones. For example, regard $G=G_2$ as a subgroup of
$SO(7)$, with its maximal torus being the subset of the maximal torus of
$SO(7)$ such that $x_1+x_2+x_3=0$. It is convenient to write
$3s_i=2x_i-x_j-x_k$ with $\{i,j,k\}=\{1,2,3\}$ and notice that
$s_1+s_2+s_3=0$. The positive roots are given by $s_i$, $1\leq i \leq 3$
and $s_j-s_k$, $j<k$, and so $\pi(s)=\sigma_3(s)\Delta(s)$. The Weyl
group $W=S_3\times\Z_2$ acts by permutations and simultaneous sign change
on the $s_i$'s. Therefore,
$L_\mu(s)=(1+(-1)^{|\mu|})\sum_{\sigma\in S_3}\sg(\sigma)(\sigma s)^\mu$
and $L_{\tau\mu}=\sg(\tau)L_\mu$ for all $\tau\in S_3$. So, taking into
account that $\sigma_1(y)=\sigma_1(s)=0$,
$$
q_\fo(s)=(\pi(y)\pi(s))^{-1}
\sum_{\lambda\in K_{m+3}}\!\frac{2\,m!}{(\lambda+\rho)!}
A_{\lambda+\rho}(y)A_{\lambda+\rho}(s)=
\sum_{\lambda\in K_{m+3}}\!\frac{2\,m!}{(\lambda+\rho)!}
\frac{S_\lambda(y)}{\sigma_3(y)}\frac{S_\lambda(s)}{\sigma_3(s)}.
$$
These invariants can be written in terms of $\sigma_2(s)\in H^4(B,\Z)$
and $\sigma_3(s^2)\in H^{12}(B,\Z)$ which form a base of the
$Ad_{G_2}$--invariant polynomials. Indeed, since $\sigma_1(s)=0$,
$\sigma_1(s^2)=-2\sigma_2(s)$ and $\sigma_2(s)^2=\sigma_2(s^2)$; see e.g.
\cite{ke}.}
\end{remark}

\section{First applications} 
In this section we prove \tref{t:tori} and \tref{t:simpleri}, and provide
several applications of the explicit expressions of the Weinstein
invariants to the case of low rank groups. In the process, we
generalize and prove some of the corollaries stated in the
Introduction.

\subsection{$G=T^n$}
In this subsection, we provide the proof of \tref{t:tori} in the
Introduction by means of a well--know algebraic result in the theory
of isometric rigidity of submanifolds.
\vspace{1.5ex}

Let $\beta:V\times V' \to W$ be a bilinear map between real vector spaces.
For $y\in V'$, define $\beta_y:V\to W$ as $\beta_y(x)=\beta(x,y)$. The set
$RE(\beta)=\{y\in V': \rk \beta_y \geq \rk \beta_z, \forall z\in V'\}$
is clearly open and dense in $V'$. The following result is essentially
contained in \cite{mo}.

\begin{lemma}\label{l:moore}
One has $\beta_z(\ker \beta_y) \subset \im \beta_y$,
for all $y\in RE(\beta), z\in V'$.
\end{lemma}
\proof
If $t$ is small, $tz+y\in RE(\beta)$. Then, $I_t=\im \beta_{tz+y}$
converges to $I_0=\im \beta_y$ as $t\to 0$. But if $x\in\ker \beta_y$,
$\beta_z(x) = \beta_{tz+y} (t^{-1}x)\in I_t$ for all $t$ small,
$t\neq 0$. Hence, $\beta_z(x) \in \im \beta_y$.
\qed
\vspace{1.5ex}

\noindent {\small \bf Proof of \tref{t:tori}:} Fix $n,s\in \N$, and let
$V\subseteq H^s(B,\R)$ be a subspace satisfying $c^k\neq 0$ for all
$c\in V\setminus\{0\}$. Consider the map
$\alpha : V \times H^{s(r-1)}(B,\R) \to H^{sr}(B,\R)$ given by
$\alpha(c,b) = c\,b$, where the product is the
cup product of the de Rham cohomology ring of $B$. We have:

\vspace{1ex}
\noindent {\it Claim. If $b\in RE(\alpha)$ and $r<k$,
the map $\alpha_b:V\to H^{sr}(B,\R)$ is a monomorphism.}
\vspace{1ex}

To prove the Claim, take $w\in \ker \alpha_b\subseteq V$.
By \lref{l:moore} we have that
$w^r = \alpha_{w^{r-1}}(w) \in \im \alpha_b$. Hence, there is $v\in V$
such that $w^r = vb$. But then $w^{r+1} = vbw = \pm v\alpha_b(w) = 0$.
Since $w\in V$, we get $w=0$ and the Claim is proved.

\vspace{1ex}
\tref{t:tori} is now a consequence of the above Claim applied to the
linear subspace $V^n\subset H^2(B,\R)$ spanned by the Chern classes of the
bundle, i.e. by the pull back of $H^2(B_{T^n})$ under the classifying map.
\qed

\vskip .3cm

\subsection{$G=U(n)$}

The general expression of the Weinstein invariants involves several
determinants and hence are difficult to use. But for certain vectors
$y\in\fg$ it can be simplified. It is thus useful to express the
nonvanishing of the invariants for some particular cases in a more
explicit way.

\begin{proposition}\label{p:2more}
Let $P$ be a $y$--fat \ $U(n)$--principal bundle over a compact
manifold $B^{2m}$ for $y\in\ft\subset\fu(n)$. Denoting
$h_k=\det(c_{j-i+1})_{1\leq i,j\leq k}$ we have:
\begin{enumerate}[\ \ \ {\rm (}a\,{\rm )}\ ]
\item If $y=(1,\dots,1)$, then $c_1^m \neq 0$;
\item If $y=(1,0,\dots,0)$, then $h_{m} \neq 0$;
\item If $y=(1+t,1,\dots,1)$, then
$\sum_{k=0}^m\binom{m+n-1}{n+k-1}t^kc_1^{m-k}h_{k} \neq 0$;
\item If $y=(1,t,0\dots,0)$, then
$\sum_{k=0}^{m/2}\binom{m+2n-3}{n+k-2}
(\sum_{i=k}^{m-k}t^i)(h_kh_{m-k}-h_{k-1}h_{m-k+1})\neq 0$.
\end{enumerate}
In particular, restrictions $(c)$ for $t=-n$ and $(d)$ for $t=-1$ also
apply for $SU(n)$--principal bundles.
\end{proposition}
\proof
Part $(a)$ follows directly from the definition of the Weinstein invariant
\eqref{e:main} since
$q_\fo(\alpha)=\int_G\la y,\alpha\ra^m dg=\int_G tr(\alpha)^m dg= c_1^m$.

For part $(b)$ apply \pref{p:main} $(a)$ to $\fo=Ad_{U(n)}(e_1)$. Using
\pref{p:si} $(g)$ and $(d)$ we obtain
$q_y=\frac{(n-1)!m!}{\rho!(m+n-1)!}h_m$ and hence $q_y \ne 0$ is
equivalent to $h_m\ne 0$.

To prove $(c)$ we use \eqref{e:main}, part $(a)$, the
proof of $(b)$ and Fubini to obtain
$$
q_y(\alpha)= \int_G\ \left( \la (e_1+\dots+e_n),\alpha\ra +
t \la e_1,\alpha\ra \right)^m dg
=\sum_{k=0}^{m} \binom{m}{k}\frac{(n-1)!k!}{\rho!(k+n-1)!}
t^k c_1^{m-k}h_k,
$$
and thus
$$
q_y(\alpha)=\frac{(n-1)!k!m!}{\rho!(m+n-1)!}\sum_{k=0}^{m}
\binom{m+n-1}{n+k-1} t^k c_1^{m-k}h_k.
$$
This proves part $(c)$.

Finally, to prove $(d)$, by \pref{p:main} $(a)$ and \pref{p:si} $(h)$, we
have $q_y=\sum_{\lambda}
\frac{m!}{(\lambda\!+\!\rho)!}S_\lambda(y)S_\lambda(x)$ for
$\lambda=(m-k,k,0,\dots,0),\ k=0,\dots,m/2$. \pref{p:si} $(a)$
and $(h)$ then imply
$$
q_y(\alpha)=\frac{(n-1)!(n-2)!m!}{\rho!(m+2n-3)!}\sum_{k=0}^{m/2}
\binom{m+2n-3}{n+k-2} t^k c_1^{m-k}h_k.
$$ which finishes our proof.
\qed

\subsection{$G=U(2)$}

We now derive a formula for the Weinstein invariants for $G=U(2)$ which
is simpler than the one obtained from \pref{p:2more} in the case of
$n=2$. This in particular proves \tref{t:simpleri} for $G=U(2)$.
\begin{proposition}\label{p:simpleru2}
Let $U(2)\to P \to B^{2m}$ be a principal bundle. If $y=(1,-1)$ is
fat then $(c_1^2-4c_2)^{m/2}\neq 0$, while if $y=(1+t,1-t)$
is fat for some $t\in\R$ then
\begin{equation}\label{e:u2eq}
\sum_{j=0}^{m/2} \binom{m+1}{2j+1}t^{2j}c_1^{m-2j}(c_1^2-4c_2)^j\neq 0.
\end{equation}
\end{proposition}
\proof
We have two fold covers
$\pi_1\colon S^1\times SU(2)\to U(2), \ (z,A)\to zA$ and $\pi_2\colon
U(2)\to S^1\times SO(3)$, obtained by dividing by $\pm Id$. All 3 have
the same polynomials $q_y$. The restrictions of
$\pi_2\circ\pi_1$ to $S^1$ and $SU(2)$ are both 2 fold covers. On the
maximal torus level we clearly have $(\pi_2\circ\pi_1)_*(a,b)=(2a,2b)$
and $(\pi_1)_*(a,b)=(a+b,a-b)$. Thus $(\pi_2)_*(s,t)=(s+t,s-t)$. Now for
the Chern class polynomials of $U(2)$ we have $c_1=s+t$ and $c_2=st$ and
thus $c_1^2-4c_2=(s-t)^2$. Hence $c_1$ becomes the Euler class for the
$S^1$ factor in $S^1\times SO(3)$ and $c_1^2-4c_2$ the Pontrjagin class
$p_1$ for $SO(3)$. For the Weinstein polynomials of $S^1$ we have
$q_y^k=y^kc_1^k$ and for $SO(3)$, using \pref{p:main}~$(b)$, $q_y^k=
\frac{2}{(k+1)} y^{k}p_1^{k/2}$ if
$k$ is even and $0$ otherwise. Since
$(\pi_2)_*(1+t,1-t)=(2,2t)$, \eqref{e:main} and Fubini imply
$$
q_y(\alpha)=\sum_{j=0}^{m/2} \binom{m}{2j}c_1^{m-2j}\frac{2^m}{2j+1} t^{2j}
(c_1^2-4c_2)^j=\frac{2^m}{m+1}\sum_{j=0}^{m/2}
\binom{m+1}{2j+1}t^{2j}c_1^{m-2j}(c_1^2-4c_2)^j.
$$
The case of $y=(1,-1)$ follows by considering $y/t=(1/t+1,1/t-1)$ and
letting $t\to\infty$ in \eqref{e:u2eq}.

As we observed in Section 1, if $P\to B$ is a $(1+t,1-t)$--fat $U(2)$
bundle, then $P/\{\pm Id\}\to B$ is a $(1,t)$--fat $S^1\times SO(3)$
bundle and the claim follows.
\qed
\vspace{.2cm}

\begin{remark}\label{r:u2}
{\rm The proof shows that $c_1^2-4c_2$ is the first Pontrjagin class of
the $SO(3)=U(2)/Z(U(2))$ bundle $P/Z(U(2))$, and $c_1$ the Euler class of
the circle bundle $P/SU(2)$.}
\end{remark}

As a consequence, we obtain the following result. The assumption is e.g.
satisfied when $b_4(B^{2m})=1$. This also proves \cref{c:cor2i} for
$U(2)$.

\begin{corollary}\label{c:cor2iu2}
Let $U(2)\to P \to B^{2m}$ be a principal bundle for which
$c_1^2=r(c_1^2-4c_2)$ for some $r\in\R$. We then have:
\begin{enumerate}[a)]
\item If $(c_1^2-4c_2)^{m/2}=0$, all Weinstein invariants vanish;
\item If $(c_1^2-4c_2)^{m/2}\neq0$ and $r=0$, there is exactly one
adjoint orbit whose Weinstein invariant vanishes;
\item If $(c_1^2-4c_2)^{m/2}\neq 0$ and $r<0$, there are exactly
$m/2$ adjoint
orbits whose Weinstein invariants vanish;
\item If $(c_1^2-4c_2)^{m/2}\neq 0$ and $r>0$, no Weinstein
invariant vanishes.
\end{enumerate}
In particular, if the bundle is fat, then $(c_1^2-4c_2)^{m/2}\neq0$
and $r>0$.
\end{corollary}
\proof
We use \pref{p:simpleru2}. For $r=0$, the vector $y=(1,1)$, i.e. $t=0$,
is clearly not fat. If $t\neq 0$, \eqref{e:u2eq} is equivalent to
$(c_1^2-4c_2)^{m/2}\neq 0$. For $r\ne 0$ we obtain
$(c_1^2-4c_2)^{m/2}\neq 0$ when $t=0$. When $t\ne 0$ we note that
$$
(1+z)^{k+1}-(1-z)^{k+1}=2z\sum_{j=0}^{[k/2]}\binom{k+1}{2j+1}z^{2j},
\ \ \ \forall z\in\C,\ k\in\N.
$$
Hence the vector $y=(1+t,1-t)$, $t\neq 0$, is fat if and only if
\begin{equation}\label{e:u2sim}
(\sqrt{r}+t)^{m+1}\neq(\sqrt{r}-t)^{m+1}.
\end{equation}
If $r>0$ this is satisfied for all $t\in\R$. If, on the contrary, $r<0$,
then one easily sees that there are exactly $m/2$ positive values of $t$
that satisfy the equality in \eqref{e:u2sim}, and thus $y=(1+t,1-t)$
cannot be fat. Notice also that $y$ and $(1-t,1+t)$ lie in the same
adjoint orbit.
\qed

\vspace{.2cm}

We point out that it is easy to state and prove similar results
to both \pref{p:simpleru2} and \cref{c:cor2iu2} for $m$ odd.

\vspace{.2cm}

For the proof of \cref{c:lens1}, let $\fs$ be the Lie algebra of
$S^1_{p,q}$, where we can assume $p\geq 1$ and $q\leq p$ by reversing the
roles of $p$ and $q$ or the orientation of the circle. Since
$\fs=\R\,(p,q)$, a straightforward computation shows that
\begin{equation}\label{e:p-q}
Ad_{U(2)}(\fs^\perp)=\R\bigcup_{t\geq|p+q|}Ad_{U(2)}(p-q+t,p-q-t).
\end{equation}

\noindent {\small \bf Proof of \cref{c:lens1}:}
If $q=p=1$, we simply get from \eqref{e:u2eq} that $(c_1^2-4c_2)^{m/2}\neq 0$.
If $q\neq p$, by \eqref{e:p-q} we need \eqref{e:u2eq} for $t\geq
|(p+q)/(p-q)|$. This is equivalent to \eqref{e:u2sim} for $t\geq
|(p+q)/(p-q)|$, which is in turn easily seen to be equivalent to
$r>-\left(\frac{1-\cos(\frac{\pi}{m+1})}{1+\cos(\frac{\pi}{m+1})}\right)
\left(\frac{p+q}{p-q}\right)^2$.
\qed

\begin{remark}\label{r:cp2}
{\rm In particular, if $m=2$, $\fs^\perp$--fatness implies that
$c_1^2=r(c_1^2-4c_2)$ with $3r>-((p+q)/(p-q))^2$. In \cite{zi1} the
results in \cite{dr} were applied to such lens space bundles as well.
It was shown that, for some orientation of the bundle, $P_+:=P/SU(2)$ is
fat and for $P_-:=P/Z(U(2)$ we have $|p_1(P_-)|<((p+q)/(p-q))^2p_1(P_+)$.
Since $p_1(P_+)=c_1^2$ and $p_1(P_-)=c_1^2-4c_2$ or vice versa,
\cref{c:lens1} for $m=2$ and the result in \cite{zi1} complement each
other. In particular, if $p+q=0$, no fat principal connection exists,
while for any other pair $p,q$ there exist fat lens space bundles over
$\CP^2$; see Section 6.}
\end{remark}

The analysis of the Weinstein invariants is easy when the dimension of
the base is small:

\vskip .2cm

{\it Fat $U(2)$--bundles over 8--dimensional manifolds}.
For $m=4$ and $y=(1+t,1-t)$ we get from \pref{p:simpleru2} $(a)$ that
$(c_1^2-4c_2)^2t^4+10c_1^2(c_1^2-4c_2)t^2+5c_1^4\neq 0$ for all $t\in \R$,
while for $y=(1,-1)$ we have $(c_1^2-4c_2)^2\neq 0$.
Therefore, no Weinstein invariant vanishes if and only if
$$
5(c_1^2(c_1^2-4c_2))^2 < c_1^4(c_1^2-4c_2)^2,\ \ \text{ or}
$$
$$
c_1^4, (c_1^2-4c_2)^2 \text{ and } c_1^2(c_1^2-4c_2)
\text{ don't vanish and have the same sign.}
$$
\vskip .2cm

\subsection{$G=SO(4)$}

As for the $U(2)$ case, for $SO(4)$ we can provide a simpler expression
for the invariants.
In particular, this proves \tref{t:simpleri} for $SO(4)$.

\begin{proposition}\label{p:simplerso4}
Let $SO(4)\to P \to B^{2m}$ be a principal bundle. If $(1,-1)$ is fat
then $(p_1+2e)^{m/2}\neq 0$, while if $(1+t,1-t)$ is fat for some
$t\in\R$, we have
$$
\sum_{j=0}^{m/2} \binom{m+2}{2j+1}t^{2j}(p_1-2e)^{m/2-j}(p_1+2e)^j\neq 0.
$$
\end{proposition}
\proof
As in the proof of \pref{p:simpleru2} we have 2-fold covers $\pi_1\colon
Sp(1)\times Sp(1)\to SO(4), \ (q_1,q_2)\to \{v\to q_1vq_2^{-1}\}$ using
multiplication of quaternions and $\pi_2\colon SO(4)\to SO(3)\times
SO(3)$, obtained by dividing by $-Id$. The restrictions of
$\pi_2\circ\pi_1$ to each $Sp(1)$ factor are again 2 fold covers and
hence $(\pi_2\circ\pi_1)_*(a,b)=(2a,2b)$. We also have
$(\pi_1)_*(a,b)=(a+b,a-b)$ and thus $(\pi_2)_*(s,t)=(s+t,s-t)$. Since
$p_1=s^2+t^2$ and $e=st$ it follows that $p_1\pm 2e = (s\pm t)^2$ are the
Pontrjagin classes of the two $SO(3)$ factors. We thus have
$$
q_y(\alpha) =2^m\sum_{k=0}^{m/2}
\binom{m}{2k}\frac{2}{m-2k+1}(p_1-2e)^{m-2k}\frac{2}{2k+1} t^{2k}
(p_1+2e))^k
$$
which, up to a factor $2^{m+2}/((m+1)(m+2))$ is the expression in
\pref{p:simplerso4}.
\qed
\vspace{.2cm}

\begin{remark}\label{r:so4}
{\rm The proof shows that $p_1\pm 2e$ is the first Pontrjagin class of
the $SO(3)=SO(4)/SU(2)_\pm$ principal bundles $P/SU(2)_\mp$, where
$SU(2)_-$ and $SU(2)_+$ are the two normal subgroups of $SO(4)$.}
\end{remark}

With the same argument as in the proof of \cref{c:cor2iu2} we easily
prove the following, which in particular generalizes \cref{c:cor2i} for
$SO(4)$.

\begin{corollary}\label{c:cor2iso4}
Let $SO(4)\to P \to B^{2m}$ be a principal bundle for which
$p_1+2e=r(p_1-2e)$ for some $r\in\R$. Then one of the following holds:
\begin{enumerate}[a)]
\item If $(p_1-2e)^{m/2}=0$, then all Weinstein invariants vanish;
\item If $(p_1-2e)^{m/2}\neq 0$ and $r=0$, then there is exactly one
adjoint orbit whose Weinstein invariant vanishes;
\item If $(p_1-2e)^{m/2}\neq 0$ and $r<0$, then there are exactly $m/2$
adjoint orbits whose Weinstein invariant vanishes;
\item If $(p_1-2e)^{m/2}\neq 0$ and $r>0$, then no Weinstein invariant
vanishes.
\end{enumerate}
In particular, if the bundle is fat, then $(p_1-2e)^{m/2}\neq0$ and $r>0$.
\end{corollary}

\begin{remark}\label{r:noso4}
{\rm In contrast to the $U(2)$ case, it is easy to see that
$\fs^\perp$--fatness implies full $\fs\fo(4)$--fatness for the
Lie algebra $\fs$ of $S^1_{p,q}\subset SO(4)$. In fact, even
$\ft^\perp$--fatness implies $\fs\fo(4)$--fatness, since
$\ft\subset Ad_{SO(4)}(\ft^\perp)$. Actually, the latter property seems
to hold for all semi--simple Lie groups of rank $>1$.}
\end{remark}

 From the dimension restriction it follows that for a fat $SO(4)$ bundle
$\dim B$ must be divisible by~$8$. In the lowest dimensional case we have:
\vspace{.2cm}

{\it Fat $SO(4)$--bundles over 8--dimensional manifolds.}
When $n=2$ and $m=4$, \pref{p:simplerso4} for $y=(1+t,1-t)$ gives
\begin{equation}\label{e:so8}
3(p_1-2e)^2+10t^2(p_1-2e)(p_1+2e)+3t^4(p_1+2e)^2\neq 0,\ \ \forall t\in\R,
\end{equation}
while for $y=(1,-1)$ we have $(p_1+2e)^2\neq 0$. Thus, no Weinstein
invariant vanishes if and only if
$$
25(p_1^2-4e^2)^2 < 9(p_1+2e)^2(p_1-2e)^2, \ \ \text{ or}
$$
$$
(p_1-2e)^2, (p_1+2e)^2 \text{ and } (p_1+2e)(p_1-2e)
\text{ don't vanish and have the same sign.}
$$

\medskip

\subsection{$G=SU(3)$}

We analyze one further case of rank $2$ groups, those with $G=SU(3)$, in
order to illustrate the difficulties one faces for other Lie groups if
one wants to express the restrictions for full fatness purely in terms of
characteristic numbers. By \eqref{e:dim}, the lowest dimensional case is
already $\dim B=32$. Here, the invariants for $SU(3)$ reduce to
\begin{eqnarray*}
c_2^2(\!\!\!\!\!\!& &\!\!\!\!\!
(511\,t^{12}+3066\,t^{11}+8814\,t^{10}+15965\,t^9+21798\,t^8+25128\,t^7+26583\,t^6+25128\,t^5\\
&+&21798\,t^4+15965\,t^3+8814\,t^2+3066\,t+511)c_2^6\\
&+&\!\!\!\!(1917\,t^{12}+11502\,t^{11}-15876\,t^{10}-184815\,t^9-498150\,t^8-757188\,t^7-834867\,t^6\\
&-&757188\,t^5-498150\,t^4-184815\,t^3-15876\,t^2+11502\,t+1917)c_3^2c_2^3\\
&+&\!\!\!\!
729(\,t^4\!+\!2\,t^3\!-\!6\,t^2\!-\!7\,t\!+\!1)(\,t^4\!+\!11\,t^3\!+\!21\,t^2\!+\!11\,t\!+\!1)(\,t^4\!-\!7\,t^3\!-\!6\,t^2\!+\!2\,t\!+\!1)
c_3^4\ )\neq 0,
\end{eqnarray*}
for all $0\leq t \leq 1$. Here we can restrict ourselves to $t\leq 1$
since if $t$ is a root, then $1/t$ also is a root. In particular,
$c_2^2(15c_3^4-21c_3^2c_2^3+c_2^6)$ and
$c_2^2(729c_3^4+1917c_3^2c_2^3+511c_2^6)$ do not vanish and have the same
sign.

Now, if $c_3^2=rc_2^3$, we write the above as $a(t)r^2+b(t)r+c(t)\neq~0$.
It is easy to see that the function $r_-(t)=(-b-\sqrt{b^2-4ac})/2a$ has
only one essential singularity in $[0,1]$ at $t_0\cong0.12920428615$, for
which $\lim_{t\to t_0^+}r_-(t)=+\infty$ and
$\lim_{t\to t_0^-}r_-(t)=-\infty$. Therefore, the function $r_-(t)$ for
$t\in[0,t_0)\cup(t_0,1]$ takes values in $(-\infty,r_2]\cup[r_1,+\infty)$,
where $r_2=r_-(0)=(-71-9\sqrt{37})/54\cong -2.3286$ and
$r_1:=r_-(1)=(15309-\sqrt{202479021})/21870$. The same argument for
$r_+(t)=(-b+\sqrt{b^2-4ac})/2a$ allows us to conclude that:
$$
\text{full fatness and $c_3^2=rc_2^3$ implies that}\
-0.30102106 \cong r_0<r<r_1\cong 0.0493593,
$$
where $r_0:=r_+(0)=(-71+9\sqrt{37})/54$.

\bigskip

A particular interesting case are $SU(3)/T^2$ fiber bundles since this is
one of the positively curved Wallach flag manifolds. But
$\ft^\perp$--fatness coincides with full $\fs\fu(3)$--fatness since it is
easy to check that $Ad_{SU(3)}(\ft^\perp)=\fs\fu(3)$. Hence we also have
that for any circle $S^1_{p,q}\subset SU(3)$ that $\fs^\perp$--fatness
implies full fatness.

\section{Fat sphere bundles} 

In this section we compute the Weinstein invariants for sphere bundles
with positive vertizontal sectional curvatures, and provide applications
related to partial fatness. We will exclude fat $S^1$--fiber bundles
which are simply in one to one correspondence with symplectic manifolds.

\subsection{Real sphere bundles}

Regard an arbitrary sphere bundle with totally
geodesic fibers of dimension $k\geq 2$ as the associated bundle to a
principal bundle $O(k\!+\!1)\to P \to B^{2m}$,
$$
\Sp^{k} \to P'=P\times_{O(k\!+\!1)} O(k\!+\!1)/O(k)\to B^{2m}.
$$
Recall that $P'$ has positive vertizontal curvatures if and only if $P$
is $\fs\fo(k)^\perp$--fat. In this situation, since
$Ad_{O(k\!+\!1)}(\fs\fo(k)^\perp)=\R Ad_{O(k\!+\!1)}(e_1)$,
\pref{p:si} together with \pref{p:main} $(b)$ and $(c)$ yield
$$
q_{\fs\fo(k)^\perp}(x)=
h_{\frac{m}{2}}(x^2)=\det(p_{j-i+1})_{1\leq i,j\leq \frac{m}{2}}\neq 0.
$$
This proves \cref{c:realspheres}. In particular, the Weinstein invariant
is independent of the dimension of the fibers and for $m\leq 8$ reduces
to:
$$
\begin{array}{lllll}
\vspace{.2cm}
\dim(B) & \ \ \ h_{\frac{m}{2}}\neq 0\\
\ \ \ \ \, 4 & p_1\neq 0\\
\ \ \ \ \, 8 & p_2\neq p_1^2\\
\ \ \ 12 & p_3\neq 2p_1p_2-p_1^3\\
\ \ \ 16 & p_4\neq p_1^4-3p_1^2p_2+2p_1p_3+p_2^2\\
\end{array}
$$
We point out that for orientable bundles, i.e., $G=SO(k+1)$, the
same formulas hold, since the term containing the Euler class vanishes.

\subsection{Complex sphere bundles}

A sphere bundle of
dimension $2n-1\geq 3$, whose underlying vector bundle has a complex
structure, can be viewed as associated to a principal bundle
$U(n)\to P \to B^{2m}$,
$$
\Sp^{2n-1} \to P'=P\times_{U(n)} U(n)/U(n\!-\!1)\to B^{2m}.
$$
Then $P'$ has positive vertizontal curvatures if and only if
$P$ is $\fu(n\!-\!1)^\perp$--fat. Since
$Ad_{U(n)}(\fu(n\!-\!1)^\perp)=\R\bigcup_{t\leq 0}Ad_{U(n)}(e_1+te_2)$,
\pref{p:2more} $(d)$ gives
\begin{equation}\label{e:nosp}
\sum_{k=0}^{[m/2]}\binom{m+2n-3}{n+k-2}\left(\sum_{i=k}^{m-k}t^i\right)
(h_kh_{m-k}-h_{k-1}h_{m-k+1})\neq 0, \ \ \forall t\leq 0,
\end{equation}
where $h_k=\det(c_{j-i+1})_{1\leq i,j\leq k}$.
In particular, for $t=0$ we obtain that
$$
\det(c_{j-i+1})_{1\leq i,j\leq m}\neq 0.
$$

For complex $\Sp^3$ fiber bundles over a 4--dimensional manifold, i.e.
$n=m=2$, fatness implies that $c_1^2(1+t+t^2)-c_2(1-t)^2\ne 0$ for all
$t\le 0$, which one easily sees is equivalent to $c_1^2=sc_2$ with $s<1$
or $s>4$. We can combine this information with the results obtained in
\cite{dr} for general $3$--sphere bundles over a 4 dimensional base. It
was shown there that there exists an orientation of the bundle such that
one of the $SO(3)$ principal bundles among $P_\pm:=P/SU(2)_\mp$, say
$P_+$, is fat and $|p_1(P_-)|<p_1(P_+)$. If the sphere bundle is a
complex sphere bundle, one has, for some choice of orientation,
$p_1(P_+)=c_1^2$ and $p_1(P_-)=c_1^2-4c_2$. The above obstruction implies
that $c_1^2=r(c_1^2-4c_2)$ with $3r>-1$ which thus complements \cite{dr}.

\medskip

\noindent {\small \bf Proof of \cref{c:zi1}:}  In
\cite{ch} and \cite{zi1} it was shown that the only $3$--sphere bundles
over $\CP^2$ that can possibly admit a fat connection
metric are the complex vector bundles with characteristic classes
$(c_1^2,c_2)=(1,1)$ or $(9,k)$, with $k=1,2,3,4$. Thus combining
both obstructions, it follows that only the sphere bundles with
$(c_1^2,c_2)= (9 ,1 )$ or $(9 ,2 )$ could possibly admit fat connection
metrics. The bundle with $(c_1^2,c_2)= (9 ,3 )$ corresponds to the
tangent bundle of $\CP^2$.
\qed

\bigskip

For $\Sp^5$ fiber bundles over 8--dimensional manifolds we have
$$
(t^4\!+\!t^3\!+\!t^2\!+\!t\!+\!1)c_1^4\!-\!3(t^4\!+\!1)c_1^2c_2\!+
\!(2\,t\!-\!1)(t\!-\!2)(1\!+\!t)^2c_1c_3\!+\!
(t^2\!-\!t\!+\!1)^2c_2^2\neq0,
$$
for all $t\le 0$, while for $\Sp^7$ fiber bundles the Weinstein
invariants are
\begin{eqnarray}\label{e:s7inv}
\notag
(t^4\!\!&+&\!\!t^3+t^2+t+1)c_1^4
+(2t^4-3t^3-13t^2-3t+2)c_1c_3
+(t-1)^4c_2^2\\
&-&\!\!(3t^2+4t+3)(t-1)^2c_1^2c_2
-(t^4-4t^3-4t^2-4t+1)c_4
\neq0,
\ \ \forall t\leq 0.
\end{eqnarray}

In particular, for $B^8=\CP^4$, the Chern classes of the tangent bundle
are $c_i=\binom{5}{i}x^i$ for a generator $x\in H^2(B,\Z)$, and thus
$(14t^4+119t^3+219t^2+119t+14)x^4\neq 0$. But this polynomial has two
real roots in $[-1,0]$, and hence $T_1\CP^4\to\CP^4$ admits no fat
connection metric. Notice that, since the sphere bundle
$T_1\CP^n\to\CP^n$ with $n\ne 1,2,4$ has no fat connection metric already
for dimension reasons, and using \cref{c:zi1}, it follows that only the
unit tangent bundle over $\CP^1$ has a fat connection metric.

\subsection{Quaternionic sphere bundles}

A sphere bundle of dimension $4n-1\geq 3$, whose underlying vector bundle
has a quaternionic structure, can be seen as an associated bundle to a
principal bundle $Sp(n)\to P \to B^{2m}$,
$$
\Sp^{4n-1} \to P'=P\times_{Sp(n)} Sp(n)/Sp(n\!-\!1)\to B^{2m}.
$$
Then, $P'$ has positive vertizontal curvatures if and only if
$P$ is $\fs\fp(n\!-\!1)^\perp$--fat. Since
$Ad_{Sp(n)}(\fs\fp(n\!-\!1)^\perp)=\R
\bigcup_{t\in\R}Ad_{Sp(n)}(e_1+te_2)$
we conclude that
$$
\sum_{k=0}^{m/4}\binom{m+4n-6}{2n+2k-3}
\left(\sum_{s=k}^{m/2-k}t^{2s}\right)
(h_{m/2-k}h_k-h_{m/2-k+1}h_{k-1})\neq 0 \ \ \forall t\in\R,
$$
where $h_k=\det(p_{j-i+1})_{1\leq i,j\leq k}$, with the $p_i$'s being the
quaternionic Pontrjagin classes. In particular for $t=0$ we get
$\det(p_{j-i+1})_{1\leq i,j\leq m/2}\neq 0$. For $n=1$, we simply obtain
$p_1^{m/2}\neq 0$ while, for $n\geq 2$, the principal bundle must be
$\fs\fp(2)\subset\fs\fp(n)$ fat and hence $32$ divides $m$.

\bigskip

The groups $G=Sp(n)\times S$ for $S=S^1,S=Sp(1)$ also act on
$\Sp^{4n-1}$ by $(A,z)\cdot v = Avz^{-1}$. Then,
$\Sp^{4n-1}=Sp(n)\times S/H$ for $H=Sp(n-1)\times \Delta S$. Thus
$Ad_G(\fh^\perp)= \R \bigcup_{0\leq t\leq 1} Ad_G(e_1-te_2,t-1)$,
and, using Fubini, we get the Weinstein invariants
$$
\sum_{i=0}^{m/2}
\sum_{k=0}^{[m/4-i/2]}\!\!\binom{m}{2i}\!\binom{m\!-\!2i\!+\!4n\!-\!6}{2n+2k-3}\!\!
\left(\!\!(t\!-\!1)^{2i}\!\!\!\!\!\sum_{s=k}^{m/2-i-k}\!\!\!\!t^{2s}\!\right)\!
w^i(h_{m/2-i-k}h_k-h_{m/2-i-k+1}h_{k-1})\neq 0
$$
for all $0\leq t\leq 1$, where $w=c_1^2$ for $S=S^1$ and $w=p_1$ for
$S=Sp(1)$.

\bigskip

\noindent {\small \bf Proof of \cref{c:s7}:}
$\Sp^7$ bundles over $\Sp^8$ are constructed by gluing two copies of
$D^8\times \Sp^7$ along the boundary $\Sp^7\times \Sp^7$ via
$(u,v)\to(u,u^kvu^l)$, where $u,v$ are unit Cayley numbers, and
$k,l\in\Z$. This defines the sphere bundle $\Sp^7\to M_{k,l}\to\Sp^8$. In
\cite{sh} it was shown that the characteristic classes of this sphere
bundle are $p_2=6(k-l)$, $e=(k+l)$. The restriction for real sphere
bundles already implies \cref{c:s7} for $k=l$. If the bundle is a sphere
bundle of a quaternionic vector bundle, we just saw that it cannot be
fat for dimension reasons. We will now determine which bundles $M_{k,l}$
carry a complex structure since any quaternionic vector bundle is also a
complex one.

 From the usual relationship between Chern and Pontrjagin classes of a
complex vector bundle it follows that $p_2=2c_4$ and $c_4=e$. Thus a
necessary condition is that $k=2l$. Admitting a complex structure is the
same as a reduction of the structure group from $SO(8)$ to $U(4)$ and
since bundles over $\Sp^8$ are classified by their gluing map along the
equator, we need to determine the image of $i_*\colon
\pi_7(U(4))\to\pi_7(SO(8))$. For this we use the long homotopy sequence
of the fibration
$$
U(4)\to SO(8)\to SO(8)/U(4)=SO(8)/SO(6)SO(2)=G_2(\R^8),
$$
where the last equality is due to one of the low dimensional isometries
of simply connected symmetric spaces. Now, $G_2(\R^8)$ is the base of
another fibration, $S^1\to V_2(\R^8)\to G_2(\R^8)$, with total space the
Stiefel manifold of $2$--frames, and they thus have the same homotopy
groups. The low dimensional homotopy groups of the Stiefel manifolds are
well known, see e.g. \cite{pa}.
In particular, $ \pi_7(G_2(\R^8))=\Z\oplus\Z_2,\
\pi_8(G_2(\R^8))=\Z_2\oplus\Z_2$ and for the homotopy groups of the Lie
groups (see e.g. \cite{mi}) we have $\pi_7(U(4))=\Z$ and
$\pi_7(SO(8))=\Z\oplus\Z$ with a basis of the latter given by the gluing
map. Thus $i_*$ is injective, and by the above its image lies in
$\Z=\{(2l,l),\ l\in\Z\}$. Since $\pi_6(U(4))=0$, the cokernel of $i_*$ is
$\Z\oplus\Z_2$ and hence $\im (i_*)= \{(2l,l),\ l \text{ even}\}$.
We conclude that the complex vector bundles are precisely the ones with
$k=2l$ for $l$ even. Notice that in this case \eqref{e:s7inv}
does not give a contradiction to fatness.

Among the complex sphere bundles, the ones that carry a quaternionic
structure are the ones for which $l$ is divisible by four. To see this,
consider the long
exact sequence in homotopy of $Sp(2)\to SU(4) \to \Sp^5$. Since
$\pi_7(\Sp^5)=\Z_2$ and $\pi_8(\Sp^5)=\Z_{24}$ (\cite{ha}), it follows
that the map from $\pi_7(Sp(2))=\Z$ to $\pi_7(SU(4))=\Z$ is
multiplication by two. This finishes our proof.
\qed

\section{Topological reduction} 

In \cite{zi1} it was conjectured that if a $G$--principal bundle admits a
fat connection, then the structure group of the bundle cannot be reduced
to any proper subgroup $H\subset G$ (where one does not assume that the
reduced bundle admits a fat connection). We make here the stronger
conjecture that this already holds when the $G$--principal connection is
$y$--fat for some $y\in\fh^\perp$. In this section we show that this is
in fact true when $H$ is a connected normal subgroup of $G$. This shows
that partial fatness, in some cases, can be used to show that the
structure groups cannot be reduced.

\bigskip

A $G$--principal bundle $\pi\colon P\to B$ is classified via its
classifying map $\phi_G\colon B\to B_G$, where $B_G=E/G$ is the
classifying
space of the Lie group $G$. The characteristic classes can then be viewed
as pull backs of cohomology classes in $H^*(B_G,\Z)$. If the structure
group of $P$ reduces, i.e., if there exists a submanifold $P'\subset P$
invariant under a subgroup $H\subset G$, then the $H$--principal bundle
$P'$ is called a {\it reduction} of $P$ and we have $P=P' \times_HG$. We
have another classifying map for $P'$, $\phi_H\colon B\to B_H$ , and
clearly $\phi_G=B_i\circ\phi_H$, where $B_i\colon B_H\to B_G$ is induced
by the inclusion map $i\colon H\to G$. Thus, if $x\in H^*(B_G,\Z)$ is a
characteristic class with $B_i^*(x)=0$, then $\phi_G^*(x)=0$ as well. In
some cases we can use the nonvanishing of certain characteristic numbers
for a fat principal bundle to show that a reduction to $H$ cannot exist.
A special case is the following result:


\begin{theorem}\label{t:red}
Let $G\to P\to B$ be a principal bundle, and $H\subset G$ a connected
normal subgroup with Lie algebra $\fh$. If the bundle reduces to $H$,
then the Weinstein invariant associated to $y$ vanishes for all
$y\in\fh^\perp$. In particular, there are no fat vectors in $\fh^\perp$.
\end{theorem}
\proof
Since $H$ is normal, $\fh$ is an ideal and thus $\fh^\perp$ is also an
ideal. Therefore, $G=H\cdot H'$ for some normal subgroup $H'\subset G$.
Since $H\times H'$ is a finite cover of $H\cdot H'$, both have the same
rational cohomology and hence the classifying spaces also have the same
cohomology. The map induced by the inclusion
$B_H=B_{H\times \{e\}}\to B_{H\times H'}=B_H\times B_{H'}$
clearly sends the characteristic classes coming
from the cohomology of $B_{H'}$ to 0. Thus the Weinstein invariant
associated to any $y\in\fh^\perp$ vanishes because of \eqref{e:peoduct}.
\qed

\begin{remark}\label{r:noso5}
{\rm If $\rk(H)=\rk(G)$, then the kernel of $B_i^*\colon H^*(B_G)\to
H^*(B_H)$ is trivial; see \cite{bo}, and thus the above method cannot be
applied. Clearly, the bigger the rank difference, the larger the kernel.
On the other hand, the bigger $\rk(G)$ is, the more difficult it is to
understand the multivariable polynomials defining the Weinstein
invariants.}
\end{remark}

\section{Examples of homogeneous fat fiber bundles} 

\medskip

Apart from the case $G=SO(2)$, where fatness is in one to one
correspondence with symplectic forms on the base $B$, the known examples
of fat bundles all arise as homogeneous bundles from inclusions $H\subset
G \subset L$, $$ G/H\to L/H\overset{\pi}{\longrightarrow} L/G=B. $$ The
metrics on $L/H$ and $L/G$ are chosen to be $L$ invariant, i.e. induced
by a left invariant metric on $L$, invariant under right translations by
$G$. If we assume that in the metric on the Lie algebra $\fl$ of $L$ the
subspaces $\fh^\perp\cap \fg$ and $ \fg^\perp \subset \fl$ (defined with
respect to a biinvariant metric) are orthogonal to each other,
Bérard--Bergery showed in \cite{bb} that the projection is a Riemannian
submersion with totally geodesic fibers. Furthermore, the submersion is
fat if and only if $[X,Y]\ne 0$ for all nonvanishing $X\in \fh^\perp\cap
\fg$ and $Y\in \fg^\perp$. In addition, Bérard--Bergery classified all
such homogeneous fat bundles.

\medskip

The above homogeneous bundle $\pi$ is associated to the $G$ principal
bundle
$$
G\to L \overset{\sigma}{\longrightarrow} L/G,
$$
since $L\times_GG/H=L/H$. If $G$ and $H$ have a normal connected subgroup
$K$ in common, and thus $G=K\cdot G'$ and $H=K\cdot H'$, we can also
choose the $G'$--principal bundle
$$
G'=G/K\to L/K \overset{\sigma}{\longrightarrow} L/G,
$$
since $L\times_{K\cdot G' }G'/H'=L/K\times_{G'}(G'/H')$. The obstructions
will be expressed in terms of the characteristic numbers $q_y$ for
$y\in\fh^\perp$ (resp. $y\in \fh'^\perp$) of the $G$ (resp. $G/K$)
principal bundle.

\bigskip

\noindent {\it Example 1: Lens Space bundles.}

\bigskip

Given the inclusion of groups
$$
U(n-1)S^1_{p,q}\subset U(n-1)U(2)\subset U(n+1),\ n\ge 2,
$$
with $S^1_{p,q}=\diag(z^p,z^q)\subset U(2)  $, we obtain the fiber bundle
over the complex Grassmannian of 2-planes in $\C^{n+1}$,
$G_2(\C^{n+1})=U(n+1)/U(n-1)U(2)$,
$$
U(n-1)U(2)/U(n-1) S^1_{p,q} \to U(n+1)/U(n-1) S^1_{p,q}\to
G_2(\C^{n+1}),
$$
with fiber $U(2)/S^1_{p,q} =SU(2)/\{\diag(z^{p},z^{q})$: $z^{p+q}=1\}$,
which is the lens space $S^3/\Z_{p+q}$ when $p+q\neq 0$. By changing the
order and replacing $z$ by $\bar{z}$ if necessary, we can assume that
$p\ge q$ and $p\ge 0$ with $\gcd(p,q)=1$. Bérard-Bergery showed in
\cite{bb} that this bundle if fat, when both the total space and the base
are equipped with a homogeneous metric, if and only if $pq>0$. We will
show now that for $pq\leq 0$ there is no fat principal connection (not
necessarily homogeneous), which will provide a proof of \cref{c:lens2}.

The above bundle can be considered to be associated to the $U(2)$
principal bundle
\begin{equation}\label{e:bundle}
U(2)\to U(n+1)/U(n-1) \to G_2(\C^{n+1}),
\end{equation}
and for the proof of \cref{c:lens2} we need its first and second Chern
classes. The cohomology ring of the base has been computed in \cite{bo}
and is given by
$$
H^*(G_2(\C^{n+1}),\Z)=
(\Z[\sigma_1,\sigma_2]\otimes Z[\bar{\sigma}_1,\dots,\bar{\sigma}_{n-1}])/
Z[\tilde{\sigma}_1,\dots,\tilde{\sigma}_{n+1}],
$$
where the $\sigma_i$'s are the symmetric polynomials in $t_1,t_2$, the
$\bar{\sigma}_i$'s the symmetric polynomials in $t_3,\dots,t_{n+1}$ and
the $\tilde{\sigma}_i$'s the symmetric polynomials in $t_1,\dots,t_{n+1}$.
Furthermore, $c_1=\bar{\sigma}_1(t_1,t_2)=t_1+t_2$ and
$c_2=\bar{\sigma}_2(t_1,t_2)=t_1t_2$ are the Chern classes of the
canonical 2--plane bundle $\xi_1$ over $G_2(\C^{n+1})$ which sends a point
into the 2--plane defining it. Similarly, $\bar{c}_1,\dots,\bar{c}_{n-1}$
are the Chern classes of the canonical $(n-1)$--plane bundle $\xi_2$ over
$G_2(\C^{n+1})$ which sends a point into the $(n-1)$--plane orthogonal
to it. Thus we can also express the cohomology ring as
$\Z[c_1,c_2]\otimes Z[\bar{c}_1,\dots,\bar{c}_{n-1}]$
divided by the relationships
$$
(1+c_1+c_2)\cup(1+\bar{c}_1+,\dots,+\bar{c}_{n-1})=
\Pi_{i=1}^{i=n-1}(1+t_i)=1,
$$
which can be regarded as the product formula for the trivial bundle
$\xi_1\oplus \xi_2$. Since
$\xi_1=U(n+1)\times_{U(n-1)U(2)}\C^2=U(n+1)/U(n-1)\times_{U(2)}\C^2$,
$c_1,c_2$ are also the Chern classes of the $U(2)$ bundle
\eqref{e:bundle}.

The above relationships imply the recursive formula
$$
\bar{c_k}=-c_1\bar{c}_{k-1}-c_2\bar{c}_{k-2},\quad k\geq1,
$$
where we set $\bar{c}_0=1$ and
$\bar{c}_{-1}=\bar{c}_{n}=\bar{c}_{n+1}=\dots=0$. Notice that this is the
same relationship as \eqref{e:inv} once we replace $\sigma_1$ by $-c_1$
and $\sigma_2$ by $c_2$. We can thus expresses the Chern classes
$\bar{c}_i$ in terms of $c_i$, as in the proof of \eqref{e:h2}, and obtain
$$
\bar{c}_k=(-1/2)^k\sum_{j=0}^{[k/2]}\binom{m+1}{2j+1}
c_1^{k-2j}(c_1^2-4c_2)^j, \qquad k\geq 1.
$$
If the lens space bundle is fat, the relationship $\bar{c}_m=0$ for
$m=2n-2$ then contradicts \eqref{e:u2eq} for $t=1$. But by \eqref{e:p-q},
\eqref{e:u2eq} is required for all $t\geq |(p+q)/(p-q)|$, and so we must
have $|(p+q)/(p-q)|>1$, or equivalently, $pq>0$, as claimed in
\cref{c:lens2}.

\bigskip

In the lowest dimensional case $n=2$, the total space is the
Aloff--Wallach space $SU(3)/S^1_{p,q}$ with embedding
$S^1_{p,q}=\diag(z^p,z^q,\bar{z}^{p+q})$, where the bundle is not only
fat, but has positive sectional curvature when $pq>0$. The metric is
obtained from the biinvariant metric on $SU(3)$ by shortening in the
direction of $U(2) =\{\diag(A,\det \bar{A}),\ A\in U(2)\}$. There are 3
such metrics corresponding to embeddings of $U(2)$ in different
coordinates and by changing the embedding, and replacing $z$ to $\bar{z}$
if necessary, any Aloff--Wallach space with $pq(p+q)\ne 0$ has a lens
space fibration with $pq>0$. From the above, it follows that for the
other two fibrations there exists no fat connection metric whatsoever. If
$n>2$, though, there is only one such fibration.

\bigskip

It is interesting to observe that, for all $n>2$, the total space admits
a metric with positive curvature on an open and dense set if $pq<0$
(\cite{wi}), and a metric with non--negative curvature and positive at
one point if $pq\ge 0$ (\cite{ta}). But these are Riemannian submersion
metrics with respect to different fibrations, where the intermediate
group $G$ in the description above is replaced by $U(n)U(1)$. They are
now fibrations over $\CP^n$ with fiber a lens space
$U(n)U(1)/U(n-1)S^1_{p,q}=U(n)/U(n-1)\cdot \Z_q=S^{2n-1}/\Z_q$.
If $q\ne 0$, there exists a metric with the above properties on the total
space, such that the projection onto $\CP^n$ is a Riemannian submersion.
But the fibers are not totally geodesic. Notice also that already from
the dimension restriction \eqref{e:dim} it follows that these bundles
cannot have a fat connection metric if $n>2$.

\bigskip

There exists another fat lens space fibration coming from the inclusions
$$
K\times S^1_{p,q}\subset K\cdot Sp(1)\times S^1\subset Q\times S^1,
$$
where $B^{4m}=Q/K\cdot Sp(1)$ is a quaternionic symmetric space and
$S^1_{p,q}\subset Sp(1)\times S^1$ is embedded with slope $(p,q)$.
Recall that a symmetric space is called {\it quaternionic} if
$Sp(1)$ acts via the Hopf action on the tangent space of the foot point.
Furthermore, in the irreducible case, each simple Lie group $Q$ gives
rise to exactly one such a space.
The above inclusions induce the fibration
$$
K\cdot Sp(1)\times S^1\!/K\times S^1_{p,q}\to Q\times S^1\!/K\times
S^1_{p,q}\to
Q\times S^1\!/K\cdot Sp(1)\times S^1=Q/K\cdot Sp(1)=B^{4m},
$$
with fiber $Sp(1)\times
S^1/\times S^1_{p,q} =Sp(1)/\{\zeta^p: \zeta\in S^1,\zeta^q=1\}$, i.e., a
lens space $S^3/Z_q$. Notice that we can assume $p\ne 0$ since otherwise
the circle acts ineffectively on the total space and the base.
Furthermore, if $q=0$, base and total space are a product with the circle
in $Sp(1)\times S^1$ and hence in both cases the bundle clearly has no
fat principal connection. Bérard--Bergery showed that this lens space
bundle is fat if and only if $pq\ne 0$.

Combining both families of examples, one sees that there exist fat lens
space bundles over $\CP^2$ for all $S^1_{p,q}$ when $p+q\ne 0$. In
\cite{zi1} it was shown that for $p+q=0$ there exists no fat connection
metric.

\bigskip
\noindent {\it Example 2: An $SO(4)$ principal bundle.}
\bigskip

Consider the $SO(4)$ principal bundle
\begin{equation}\label{e:G2}
SO(4)\to G_2\to G_2/SO(4).
\end{equation}
Bérard--Bergery showed that this bundle is
$\fsu(2)_\pm^\perp$ fat where $SU(2)_\pm$ are the two normal subgroups of
$SO(4)$ corresponding to the image of $S^3\times\{e\}$ and $\{e\}\times
S^3$ under the two fold cover $S^3\times S^3\to SO(4)$. Thus the
associated bundles $ G_2\times_{SO(4)}SO(4)/SU(2)_\pm=G_2/SU(2)_\pm\to
G_2/SO(4)$ are both fat $SO(3)=SO(4)/SU(2)_\pm$ principal bundles.

We now compute $p_1$ and $e$ of the $SO(4)$ bundle \eqref{e:G2}.
For this we use the Borel method which we now recall. We have a
commutative diagram

\begin{picture}(100,160)\label{dddd2}  
\put(145,130){$G$}
\put(258,130){$G$}
\put(180,133){\vector(1,0){55} }
\put(145,117){ \vector(0,-1){20} }
\put(258,117){  \vector(0,-1){20} }
\put(140,80){$ G/H$ }
\put(258,80){$E$}
\put(180,83){\vector(1,0){55} }
\put(145,65){  \vector(0,-1){20} }
\put(150,55){ {$\varphi_{H}$ } }
\put(258,65){  \vector(0,-1){20} }
\put(145,30){$B_{H}$}
\put(257,30){$B_{G}$}
\put(180,33){\vector(1,0){55} }
\put(200,22){$B_{i}$}
\end{picture}  

\noindent where $i$ is the inclusion $i \colon H\to G$. Thus the left
hand side $G$ principal bundle is the pull back of the universal bundle
on the right. The differentials in the universal spectral sequence are
well known and the ones in the left hand side fibration are induced by
naturality as soon as we know the map in cohomology
$B_i^*\colon H^*(B_G)\to H^*(B_H)$. The map $\varphi_{H}$ is the
classifying map of the $H$ principal bundle $H\to G\to G/H$ which can be
determined by the edge homomorphism in the spectral sequence and this
will then give us the values of the characteristic classes of the $H$
principal bundle.

In order to compute $B_i^*$, we let $T_G\subset G$ and
$T_H\subset H$ be maximal tori and use the commutativity of the diagram

\begin{picture}(100,105) 
\put(153,75){$B_H$}
\put(240,75){$B_G$}
\put(180,78){\vector(1,0){55}}
\put(200,67){$B_i$}
\put(155,60){\vector(0,-1){20}}
\put(160,50){$B_{j_H}$}
\put(240,60){\vector(0,-1){20}}
\put(250,50){$B_{j_G}$}
\put(153,25){$B_{T_H}$}
\put(240,25){$B_{T_G}$}
\put(180,28){\vector(1,0){55}}
\put(200,17){$B_i$}
\end{picture}

\noindent
We choose coordinates $(t_1,\dots,t_n)$ of the (integral lattice of the)
maximal torus $T_G\subset G$ and, by abuse of notation, let $t_i\in
H^1(T_{G},\Z)=Hom(\pi_1(G),\Z)$ and hence $\bar{t}_i\in H^2(B_{T_{G}})$
via transgression in the spectral sequence of the universal bundle of
$T_G$. We then have $H^*(B_{T_{G}})=P[\bar{t}_1\dots , \bar{t}_n]$ and
$B_{j_G}^*$ is injective on the torsion free part of $H^*(B_G)$ with
image $H^*(B_{T_{G}})^{W_{G}}$, where $W_{G}$ is the Weyl group of $G$,
and similarly for $H$. We thus only need to compute $B_i^*\colon
H^*(T_G)\to H^*(T_H)$, which is easily done.

We now apply all this to $G=G_2$ and $H=SO(4)$. Additional complications
arise since the cohomology of the groups and their classifying spaces
contain torsion. In \cite{fe2} it was shown that for both $SO(4)$ and
$G_2$, the homomorphism $B_{i_G}^*\colon H^*(B_G,\Z)\to H^*(B_{T_G},\Z)$,
after dividing by the torsion groups, is injective with image the Weyl
group invariant subalgebra. For $G=SO(4)$, if we use the coordinates for
$T_G$ as in Section 1, the transgression $\bar{x}_i\in
H^2(B_T,\Z)=\Z\oplus\Z$ form a basis, and $p_1=\bar{x}_1^2+\bar{x}_2^2$
and $e=\bar{x}_1\bar{x}_2$ form a basis of the Weyl group invariant
subalgebra. Here $p_1$ and $e$ are the universal Pontrjagin and Euler
classes. There are elements $y_1,y_2\in H^3(G,\Z)=\Z\oplus\Z$ such that
$d_3(y_1)=p_1$ and $d_3(y_2)=e$ in the universal spectral sequence for
$SO(4)$.

The maximal torus of $G_2$ is given by $(t_1,t_2,t_3)$ with $\sum t_i=0$
and a basis of the Weyl group invariant algebra is $x=\frac12
\sigma_1(\bar{s}_i^2)$ and $y=\sigma_3(\bar{s}_i^2)$ where
$s_i=\frac13(2t_i-t_j-t_k)$. Since the positive roots are $s_i$, $1\leq i
\leq 3$ and $s_j-s_k=t_j-t_k$, $j<k$, one easily sees that the roots
$s_3,\ s_1-s_2$ span a subalgebra isomorphic to the Lie algebra of the
(unique) $SO(4)$ in $G_2$. In terms of $x_1,\ x_2$, the roots are $\pm
x_1\pm x_2$ and hence we can choose $x_1=s_1,\ x_2=s_2$. Thus
$p_1=s_1^2+s_2^2,\ e=s_1s_2$ and since $x=\frac12
\sigma_1(\bar{s}_i^2)=s_1^2+s_2^2+s_1s_2$ it follows that $d_4(x)=p_1+e$
in the spectral sequence of the left hand side fibration in the first
diagram. Thus $H^4(G_2/SO(4),\Z)=\Z[a]$ with $e=a,\ p_1=-a$ and by
Poincaré duality $H^8(G_2/SO(4),\Z)=\Z[a^2]$. Thus $p_1+2e=r(p_1-2e)$
with $r=-1/3$. Hence the $SO(4)$ principal bundle cannot be fat. Using
\eqref{e:so8}, it follows that the Weinstein invariant for $y=(1+t,1-t)$
is
$$
(9-10t^2+t^4)(p_1-2e)^2\neq 0,
$$
with zeroes $t=\pm 1, \pm 3$ and hence the Weinstein invariant is $0$
for the two adjoint orbits with $y=(1,0)$ and $y=(1,-2)$. Notice that the
fat bundle by Bérard--Bergery has $y=(1,\pm 1)$. Thus we have
\begin{corollary}\label{c:g2}
The principal bundle $SO(4)\to G_2\to G_2/SO(4)$ has a homogeneous
connection metric which is $(1,\pm 1)$--fat, but admits no
$y$--fat principal connection for $y=(1,0)$ or $y=(1,-2)$. In particular,
the associated 3--sphere bundle $\Sp^3\to G_2/SO(3)\to G_2/SO(4)$ does not
admit a fat connection metric.
\end{corollary}

\bigskip
\noindent {\it Example 3: Sphere Bundles.}
\bigskip

All remaining examples of fat homogeneous fibrations in \cite{bb} are
bundles with fiber $\Sp^n$ or $\R P^n$ represented as $SO(n)/SO(n-1)$ or
$SO(n)/O(n-1)$. We describe next a typical case.

\medskip

The inclusions $Sp(1)Sp(1)Sp(n-2)\subset Sp(2)Sp(n-2)\subset Sp(n)$
induce the fibration
$$
\Sp^4\to M\to G_2(\QH^n).
$$
The $Sp(2)$ principal bundle $ Sp(2)\to Sp(n-2)/Sp(2) \to G_2(\QH^n)$ has
quaternionic Pontrjagin classes $p_1,p_2$ with $H^*(G_2(\QH^n),\Z)=
\Z[p_1 ,p_2]/\{(1+p_1+p_2) \sum_{i=0}^{i=n-2}\bar{p}_i =1\}$.
If $n=3$ the total space is a
positive curved Wallach flag manifold with base $\HP^2$. Thus
$p_1^2=p_2=1$ in this special case. But as an $\Sp^4$ bundle the
structure group is $SO(5)$ and the 2-fold cover $Sp(2)\to SO(5)$ induces
a map $H^*(B_{SO(5)})\to H^*(B_{Sp(2)})$ which relates the real and
quaternionic Pontrjagin classes. Using this, one easily shows that the
real Pontrjagin numbers are $p_1^2=4$ and $p_2=-3$, which is consistent
with the obstruction for fat real sphere bundles described above.

{\renewcommand{\baselinestretch}{1}
\hspace*{-25ex}
\begin{tabbing}
\indent \= IMPA -- Estrada Dona Castorina, 110\\
\> 22460-320 -- Rio de Janeiro -- Brazil\\
\> email: luis@impa.br\\
\end{tabbing}}

\vspace*{-7ex}

{\renewcommand{\baselinestretch}{1}
\hspace*{-25ex}
\begin{tabbing}
\indent \= University of Pennsylvania\\
\> Philadelphia, PA 19104, USA\\
\> email: wziller@math.upenn.edu\\
\end{tabbing}}

\vspace*{-7ex}

\end{document}